\documentclass[11pt,twoside]{amsart}

\usepackage[dvipsnames]{xcolor}
\usepackage{hyperref}
\hypersetup{colorlinks=true, linkcolor=Blue, citecolor=Maroon}
\usepackage{amsthm, amsmath, amscd, amssymb,centernot,txfonts}
\usepackage{tikz-cd}

\usepackage{lipsum}
\usepackage{multicol}

\usepackage[normalem]{ulem}
\usepackage{amsfonts}
\usepackage{mathrsfs}
\usepackage{parskip}
\usepackage[all]{xy}
\usepackage[left=2cm,top=2.5cm,bottom=3cm,right=3cm]{geometry}
\usepackage{dirtytalk}
\usepackage{verbatim}
\usepackage{mathtools}
\usetikzlibrary{shapes,arrows}
\usepackage{etoolbox}
\usepackage{dynkin-diagrams}

\usepackage{enumitem}
\usepackage{chngpage}

\usepackage{adjustbox}
\usepackage[utf8]{inputenc}
\usepackage{fourier} 
\usepackage{array}
\usepackage{makecell}

\usepackage[linktocpage=true]{}

\usepackage{color, colortbl}
\definecolor{LightCyan}{rgb}{0.88,1,1}

\usepackage[normalem]{ulem}

\usepackage[first=0,last=9]{lcg}

\newcommand{\floor}[1]{\left\lfloor #1 \right\rfloor}

\newcommand{\stkout}[1]{\ifmmode\text{\sout{\ensuremath{#1}}}\else\sout{#1}\fi}

\theoremstyle{plain}
\numberwithin{equation}{section}
\newtheorem{theorem}{Theorem}[section]
\newtheorem{proposition}[theorem]{Proposition}
\newtheorem{lemma}[theorem]{Lemma}

\newtheorem{set-up}[theorem]{Set-up}

\theoremstyle{definition}
\newtheorem{remark}[theorem]{Remark}

\newtheorem{example}[theorem]{Example}
\newtheorem{definition}[theorem]{Definition}

\newcommand*{\QEDA}{\hfill\ensuremath{\blacksquare}}
\newcommand*{\QEDB}{\hfill\ensuremath{\square}}

\tikzstyle{decision} = [diamond, draw, , 
    text width=4.5em, text badly centered, node distance=3cm, inner sep=0pt]
\tikzstyle{block} = [rectangle, draw, , 
    text width=10em, text centered, rounded corners, minimum height=2em]
\tikzstyle{block1} = [rectangle, draw, , 
    text width=5em, text centered, rounded corners, minimum height=2em]    
\tikzstyle{line} = [draw, -latex']
\tikzstyle{cloud} = [draw, ellipse,, node distance=3cm,
    minimum height=2em]

\begin{document}
\title[Smoothing of multiple structures on embedded Enriques manifolds]{Smoothing of multiple structures on embedded Enriques manifolds}

\author[J. Mukherjee]{Jayan Mukherjee}
\address{Department of Mathematics, University of Kansas, Lawrence, KS 66045}
\email{j899m889@ku.edu}

\author[D. Raychaudhury]{Debaditya Raychaudhury}
\address{Department of Mathematics, University of Kansas, Lawrence, KS 66045}
\email{debaditya@ku.edu}

\subjclass[2020]{14C05, 14D06, 14D15, 14E20, 53C26}
\keywords{Ropes, Enriques manifolds, hyperk\"ahler manifolds, Calabi-Yau manifolds, smoothing of multiple structures, Hilbert scheme}

\maketitle
\begin{abstract}
We show that given an embedding of an Enriques manifold of index $d$ in a large enough projective space, there will exist embedded multiple structures with conormal bundle isomorphic to the trace zero module of the universal covering map, the universal cover being either a hyperk\"ahler or a Calabi-Yau manifold. We then show that these multiple structures (also known as $d$-ropes) can be smoothed to smooth hyperk\"ahler or Calabi-Yau manifolds respectively. Hence we obtain a flat family of hyperk\"ahler  (or Calabi-Yau) manifolds embedded in the same projective space which degenerates to an embedded $d$-rope structure on the given Enriques manifold of index $d$. The above shows that these $d$-rope structures on the embedded Enriques manifold are points of the Hilbert scheme containing the fibres of the above family. We show that they are smooth points of the Hilbert scheme when $d=2$.
\end{abstract}
\section{Introduction} In this article we deal with two topics namely, deformations of cyclic \'etale covers of some non-simply connected projective  varieties embedded in a projective space and multiple structures (also known as ropes) supported on the same variety inside the same projective space.\par 
A multiple structure or more precisely a $d$-rope on a reduced algebraic variety $X$ is an everywhere non-reduced scheme $Y$ supported on $X$, such that the ideal sheaf $\mathscr{I}$ of the embedding of $X$ in $Y$ is of square zero and $\mathscr{I}$ is a locally free $\mathscr{O}_X$  module of rank $d$. Ropes of rank two (also known as ribbons) first appeared in the context of Green's canonical syzygy conjecture (see \cite{Green}). Bayer and Eisenbud (see \cite{BE}) defined ribbons over schemes, as a result of which they formulated the canonical syzygy conjecture for ribbons (see \cite{Deo}) whose solution implies the canonical syzygy conjecture in the generic case. The double structures as defined by Bayer and Eisenbud often appear as the flat limits of smooth algebraic varieties. More precisely we call a $d$-rope $D$ \textit{smoothable}  if there exists a flat family of schemes over a curve $T$ such that the central fibre is $D$ and after possibly shrinking $T$ we can assume all other fibres are smooth. Results on smoothings of ribbons over curves can be found in \cite{Fong} and \cite{Gon}. Smoothings of $2$-ropes on surfaces, also called \textit{carpets} were studied by Gallego, Gonz\'alez and Purnaprajna in \cite{GGP} and \cite{GP}. In the above papers they prove strong smoothing results on certain special carpets called $K3$ carpets on Enriques surfaces which are regular and have trivial canonical bundle. The above authors in the series of papers \cite{GGP2}, \cite{GGP3}, \cite{GGP4} further show the connection between existence and smoothings of multiple structures of conormal bundle $\mathscr{E}$ and deformation of finite covers with trace zero module $\mathscr{E}$. In \cite{GGP1} the authors give an infinitesimal condition to ensure that a finite morphism can be  deformed to an embedding. The criterion is entirely cohomological and demands that the tangent vector of a deformation along a curve correspond to a surjective homomorphism between the conormal bundle and the trace zero module. \par
At this point we discuss \'etale covers of some non-simply connected projective algebraic varieties. Enriques surfaces are a class of surfaces that admit an \'etale  double covering (which is in fact universal) by $K3$ surfaces. They are minimal surfaces of Kodaira dimension zero with $p_g = 0$ and $q = 0$. In \cite{GGP}, the authors  obtain a smoothing family of $K3$ surfaces of a so-called ``$K3$ carpet'' on an Enriques surface, thereby showing that such carpets arise naturally as flat limits of $K3$ surfaces. The geometry of the analogue of $K3$ surfaces in higher dimensions, which are known as hyperk\"ahler and Calabi-Yau manifolds, are central to the study of higher-dimensional varieties. In 2011 Boissi\`ere, Nieper and Sarti and in the same year Oguiso and Schr\"oer independently constructed higher-dimensional analogues of Enriques surfaces (see \cite{BNS}, \cite{OS}). Although their definitions differ a little, in this article the definition of an Enriques manifold that we take is either an Enriques manifold in the sense of Oguiso and Schr{\"o}er or in the sense of Boissi\`ere, Nieper and Sarti. In any case, such higher-dimensional Enriques manifolds are smooth projective non-simply connected (as a complex manifold) varieties in even dimension with torsion canonical bundle whose universal covers are either hyperk\"ahler or Calabi-Yau manifolds. The varieties were obtained by finding fixed point free actions on either hyperk\"ahler or Calabi-Yau manifolds and the existing examples have indices either $2$, $3$ or $4$ (index of an Enriques manifold is the degree of its universal covering morphism). Note that if the universal cover $X$ is a Calabi-Yau manifold then the only index possible is $2$ since $\chi(\mathscr{O}_X) = 2$. The latter authors in \cite{OS1} also describe the period map and prove the local Torelli theorem for higher-dimensional Enriques manifolds. \par 
In this paper we investigate for an Enriques manifold $Y$ of index $d$ constructed in either of \cite{BNS} or \cite{OS}, the deformations of the composite of the universal covering map of $Y$ with an embedding $i: Y \hookrightarrow \mathbb{P}^N$ with the aid of the infinitesimal condition described in \cite{GGP1}. 
For an Enriques manifold $Y$ with canonical bundle $K_Y$ which is torsion of index $d$, 
$\mathscr{E}$ will always denote the following vector bundle,
\begin{align*}
    \mathscr{E}=\bigoplus\limits_{i=1}^{d-1}K_Y^{\otimes i}.
\end{align*}
We show that for an embedding of an Enriques manifold $Y$ of index $d$ and dimension $2n$ inside a projective space of dimension $N$, the composite $\varphi$ of the covering map and the embedding can be deformed to a finite birational morphism if $N \geq 2n+\floor{d/2}$ and $\mathscr{O}_Y(1) \otimes K_Y^{\otimes i}$ is base point free for all $1\leq i\leq d-1$. Under the same global generation condition, if $N \geq 4n+d-1$, we show that $\varphi$ can be deformed to an embedding by showing the existence and smoothing of embedded $d$-rope structures supported on $Y$ with conormal bundle $\mathscr{E}$. Hence in the latter case, we obtain a smoothing family as the image of a flat family of morphisms from hyperk\"ahler or Calabi-Yau manifolds to a projective space whose general fibre is an embedding. This shows how these embedded $d$-rope structures on Enriques manifolds can be naturally thought of as degenerations of families of hyperk\"ahler or Calabi-Yau manifolds inside a projective space. In the following paragraph we briefly summarize the contents of the article. \par 
Section ~\ref{prelims} deals with the notation and preliminaries. In section ~\ref{RE}, we study embedded $d$-ropes with conormal bundle $\mathscr{E}$ on an Enriques manifold of dimension $2n$ inside a projective space of dimension $N$ in detail, and compute the dimension of the quasi-projective space parametrizing such $d$-ropes (see Propositions ~\ref{countinghk}, ~\ref{countingcy}). This section hence answers the question of existence of embedded $d$-ropes with conormal bundle $\mathscr{E}$. In Section ~\ref{smoothing}, 
we apply the results of Gallego, Gonz\'alez and Purnaprajna and prove our main results on smoothing of such ropes (see Theorem ~\ref{smoothinghk}, Theorem ~\ref{smoothingcy}). Since we have an embedded smoothing, we see that such ropes on Enriques manifolds correspond to points in the Hilbert scheme of embedded hyperk\"ahler or Calabi-Yau manifolds and hence a natural question is the smoothness of the points. We answer this question at the end of Section ~\ref{smoothing} (see Theorem ~\ref{smhs}) for ropes of index $d = 2$ (we call them  $K$- trivial ribbons). In the Appendix (Section ~\ref{App}) we prove that if an Enriques manifold $Y$ of index $d = 2$ whose universal cover is one of the known examples of hyperk\"ahler six-folds is embedded inside $\mathbb{P}^N$, then $N\geq 13$.
\subsection*{Acknowledgements.} We are grateful to our advisor Professor B. P. Purnaprajna for suggesting to us this problem and for supporting us throughout this work. We thank Professor Keiji Oguiso for his comments on an earlier draft of this paper and for pointing out to us the work of Taro Hayashi (\cite{Ha}). We thank the anonymous referee and Professor Olivier Debarre for helpful suggestions and corrections that substantially helped in the improvement of the exposition.

\section{Preliminaries and Notation} \phantomsection\label{prelims}Throughout this paper we shall work over the field of complex numbers $\mathbb{C}$. In this section, first we will recall some basic results on Calabi-Yau, hyperk\"ahler and Enriques manifolds. Then we will define ropes and state their characterizations. Ribbons are special kinds of ropes and we will encounter regular $K$- trivial ribbons that can be thought of as a higher-dimensional analogue of $K3$ carpets. 
\subsection{Calabi-Yau and hyperk\"ahler manifolds.} \phantomsection\label{prelimcyhk}A $K3$ surface $S$ is by definition a compact, complex, K\"ahler manifold which is regular and has trivial canonical bundle. Note that the definition is equivalent to having a holomorphic symplectic form on $S$. However in higher dimensions these two notions do not coincide giving rise to Calabi-Yau and hyperk\"ahler manifolds. \par
\begin{definition}
 \phantomsection\label{defcy} A compact K\"ahler manifold $M$ of dimension $n\geq 3$ is called {\it Calabi-Yau} if it has trivial canonical bundle and the Hodge numbers $h^{p,0}(M)$ vanish for all $0<p<n$.
\end{definition} 
With this definition, Calabi-Yau manifolds are necessarily projective but the following definition of hyperk\"ahler manifolds does not imply projectivity. There are examples of non-projective hyperk\"ahler manifolds (see for example \cite{Gr}).
\begin{definition}
\phantomsection\label{defhk}A compact K\"ahler manifold $M$ is called {\it hyperk\"ahler} if it is simply connected and its space of global holomorphic two forms is spanned by a symplectic form. 
\end{definition} 
By the decomposition theorem of Bogomolov (see \cite{Bo}), any complex manifold with trivial first Chern class admits a finite \'etale cover isomorphic to a product of complex tori, Calabi-Yau manifolds and hyperk\"ahler manifolds. These spaces can be thought of as the ``building blocks" for manifolds with trivial first Chern class. \par 
Let $X$ be a hyperk\"ahler manifold. The symplectic form $\sigma$ induces an isomorphism $T_X\to \Omega_X$ by sending $\delta$ to $\delta'\to \sigma(\delta,\delta')$. Moreover, the symplectic form ensures that $K_X$ is trivial and $\textrm{dim}(X)$ is even. It is also known that $\chi(\mathscr{O}_X)=n+1$ and the following values of $h^p(X,\mathscr{O}_X)$
\[h^p(X,\mathscr{O}_X)=
    \begin{dcases}
        1 & \textrm{$p$ is even,} \\
        0 & \textrm{$p$ is odd.} \\
    \end{dcases}
\]
There are many families of examples for Calabi-Yau varieties but only a few classes of examples for hyperk\"ahler varieties are known. Beauville first produced examples of two distinct deformation classes of compact hyperk\"ahler manifolds in all even dimensions greater than or equal to $2$ (see \cite{Be}). The first example is the Hilbert scheme $S^{[n]}$ of length $n$ subschemes on a $K3$ surface $S$. The second one is the generalized Kummer variety $K^n(A)$ which is the fibre over the $0$ of an Abelian surface $A$ under the morphism $\phi\circ\psi$ (see the diagram below)
\begin{align*}
    A^{[n+1]} \xrightarrow{\psi} A^{(n+1)}\xrightarrow{\phi} A
\end{align*}
where $A^{[n+1]}$ Hilbert scheme of length $n+1$ subschemes on the Abelian surface $A$, $A^{(n+1)}$ is the symmetric product, $\psi$ is the Hilbert-Chow morphism and $\phi$ is the addition on $A$. Two other distinct deformation classes of hyperk\"ahler manifolds are given by O'Grady in dimensions $6$ and $10$ which appear as desingularizations of certain moduli spaces of sheaves over symplectic surfaces (see \cite{OG1}, \cite{OG2}). All other known examples are deformation equivalent to one of these.\par 
We recall the definitions of the Beauville form and the Fujiki constant on hyperk\"ahler manifolds (see \cite{Be} and \cite{Fu}). 

\begin{theorem}
\phantomsection\label{BFFC}Let $X$ be a hyperk\"ahler variety of dimension $2n$. There exists a quadratic form $q_X:H^2(X,\mathbb{C})\rightarrow \mathbb{C}$ and a positive constant $c_X\in\mathbb{Q}_+$ such that for all $\alpha$ in $H^2(X,\mathbb{C})$, $\int_X\alpha^{2n}=c_X\cdot q_X(\alpha)^n$. The above equation determines $c_X$ and $q_X$ uniquely if one assumes the following two conditions.
\begin{itemize}
    \item[(1)] $q_X$ is a primitive integral quadratic form on $H^2(X,\mathbb{Z})$,
    \item[(2)] $q_X(\sigma,\bar{\sigma})>0$ for all $0\neq\sigma\in H^{2,0}(X)$.
\end{itemize}
Here $q_X$ and $c_X$ are called the Beauville form and the Fujiki constant respectively.
\end{theorem}
The Beauville form and the Fujiki constant are the fundamental invariants of a hyperk\"ahler variety. Moreover, they help finding an explicit form of the Riemann-Roch theorem for hyperk\"ahler varieties.\par
 \begin{theorem}
 \phantomsection\label{GRR} (See \cite{Fu}, \cite{Gr}) Let $X$ be a hyperk\"ahler variety of dimension $2n$. Assume that $\alpha\in H^{4j}(X,\mathbb{C})$ of type $(2j,2j)$ on all small deformations of $X$. Then there exists a constant $C(\alpha)\in \mathbb{C}$ depending on $\alpha$ such that $\int_X\alpha\cdot\beta^{2n-2j}=C(\alpha)\cdot q_X(\beta)^{n-j}$ for all $\beta\in H^2(X,\mathbb{C})$.
 \end{theorem}
As a consequence of the theorem above, we get the following form of the Riemann-Roch formula for a line bundle $L$ on a hyperk\"ahler variety of dimension $2n$ (see \cite{Hu}),
\begin{align}\phantomsection\label{SRR}
    \chi(X,L)=\sum_{i=0}^{n}\dfrac{a_i}{(2i)!}q_X(c_1(L))^i
\end{align} where $a_i=C(\textrm{td}_{2n-2i}(X))$. Here $a_i$'s
 are constants depending only on the topology of $X$ (see for example the discussion after Theorem 3.1 in \cite{CaoJ}). We refer to \cite{Hu} for an introduction to the theory of hyperk\"ahler manifolds.\par 

\subsection{Enriques manifolds.}\phantomsection\label{prelimenriques}
Enriques surfaces are defined as non-simply connected surfaces $S$ who-\\se universal cover is a $K3$ surface. Generalization of Enriques surfaces to higher dimensions has been constructed by Boissi\`ere, Nieper-Wisskirchen, Oguiso, Sarti and Schr\"oer (see \cite{BNS}, \cite{OS}). We take the following definition of Enriques manifolds.
\begin{definition}\phantomsection\label{defenriques}
Let $Y$ be a connected, compact, complex, smooth, K\"ahler manifold of dimension $2n$. $Y$ is called an \textit{Enriques manifold} if there exists an integer $d\geq 2$, called the index of $Y$, such that the canonical bundle $K_Y$ has order $d$ in the Picard group $\textrm{Pic}(Y)$ of $Y$, the fundamental
group $\pi(Y )$ is cyclic of order $d$ and the holonomy representation of its universal cover is irreducible.
\end{definition}
Note that when $d>2$ the manifold $Y$ is an Enriques manifold in the sense of \cite{OS}. When $d=2$ the manifold $Y$ is either an Enriques variety in the sense of \cite{BNS} if its universal cover is a Calabi-Yau manifold or it is an Enriques manifold in the sense of \cite{OS} if its universal cover is a hyperk\"ahler manifold. Consequently, $Y$ is projective and so is its universal cover (for the proof of this statement see \cite{BNS} Proposition 2.1 and \cite{OS}, Corollary 2.7). \par 
\indent Let $Y$ be an Enriques manifold of index $d$ and let $\pi:X\to Y$ be its universal cover. Then, 
\begin{align*}
    \pi_*\mathscr{O}_X=\mathscr{O}_Y\oplus K_Y\oplus K_Y^{\otimes 2}\oplus\cdots\oplus K_Y^{\otimes d-1}=\mathscr{O}_Y\oplus \mathscr{E}.
\end{align*}
Here and in the subsequent sections, $\mathscr{E}$ will always stand for $K_Y\oplus K_Y^{\otimes 2}\oplus\cdots\oplus K_Y^{\otimes d-1}$. If $X$ is a Calabi-Yau manifold, $h^i(\mathscr{O}_Y)=0$ for $1\leq i\leq \textrm{dim}(Y)-1$, $\chi(\mathscr{O}_Y)=1$ and $d=2$. If $X$ is hyperk\"ahler then $\chi(\mathscr{O}_Y)=(\textrm{dim}(Y)+2)/2d$ and the following gives the cohomology of the structure sheaf. $h^p(Y,\mathscr{O}_Y)$,
\[h^p(Y,\mathscr{O}_Y)=
    \begin{dcases}
        1 & 2d\mid p \textrm{ and } p\leq \textrm{dim}(Y), \\
        0 & \textrm{otherwise.} \\
    \end{dcases}
\]
\indent The examples below are taken from \cite{OS}, they provide examples of Enriques manifolds of index two. In general, it is known that there are Enriques manifolds of index 2, 3, 4 (see \cite{BNS}, Theorem 2.2 and \cite{OS}, Theorem 6.5). The first one is an Enriques manifold whose universal cover is a Calabi-Yau manifold. The second and third examples are Enriques manifolds whose universal cover is a hyperk\"ahler manifold of the form $K3^{[n]}$ and $K^n(A)$ respectively.
\begin{example}\phantomsection\label{ex1}
Let $S$ be an Enriques surface and $S^{[n]}$ be the Hilbert scheme of $n$ points on $S$. Then $S^{[n]}$ is an Enriques manifold of index two whose universal cover is a Calabi-Yau manifold (see \cite{OS}, Theorem 3.1).
\end{example}
\begin{example}\phantomsection\label{ex2}
Let $S'$ be an Enriques surface and $S\to S'$ be its universal covering. Then $S$ is a $K3$ surface with a free $G=\mathbb{Z}/2\mathbb{Z}$ action on it. Suppose $n$ is odd. Then the induced action of $G$ on the Hilbert scheme $S^{[n]}$ is free and $Y:=S^{[n]}/G$ is an Enriques manifold of index two whose universal cover is the hyperk\"ahler manifold $S^{[n]}$ (see \cite{OS}, Proposition 4.1).
\end{example}
\begin{example}\phantomsection\label{ex3}
Let $E$ and $E'$ be two elliptic curves, $a'\in E'$ be a point of order two and $a\in E$ be an arbitrary point. Then $(b,b')\mapsto (-b+a,b'+a')$ is a free involution on $A:=E\times E'$ which induces a free $G:=\mathbb{Z}/2\mathbb{Z}$ action on the Hilbert scheme $A^{[n]}$. Then $G$ leaves $K^n(A)$ invariant and the induced action on $K^n(A)$ is free. Moreover, $K^n(A)/G$ is an Enriques manifold of index two whose universal cover is $K^n(A)$ (\cite{OS}, Proposition 4.2).
\end{example}
\subsection{Deformations of the universal covering morphism} \phantomsection\label{defsp} This subsection provides background and motivation for the study that we do in the next subsection. Let $Y$ be an Enriques manifold. Let $\pi:X\to Y$ be the covering map from the universal cover $X$. 
Let us address the following question:\par 

\noindent \textit{Can a general small deformation of $X$ be realized again as an universal cover of a deformation of $Y$?}\par  

It is interesting to note the dichotomy of the answer to the above question depending whether the universal cover is hyperk\"ahler or Calabi-Yau. It is clear from \cite{OS} and \cite{OS1} that this is never the case when $X$ is a hyperk\"ahler manifold. However it has been shown by Hayashi in \cite{Ha} that if the universal cover of the Enriques manifold is Calabi-Yau, then the above question has a positive answer. We include  proofs of these statements in Remark ~\ref{defsphk} and Theorem ~\ref{defspcy} for the sake of completeness. The proof of Theorem ~\ref{defspcy} was independently done by us and differs from the original proof of \cite{Ha} by the use of \cite{Weh}, Corollary $1.11$.
\begin{remark}\phantomsection\label{defsphk}
 Let $Y$ be an Enriques manifold of degree $d$. Let $X \xrightarrow{\pi} Y$ be the universal cover of $Y$ where $X$ is a hyperk\"ahler manifold. Let $\mathscr{B}_X$ and $\mathscr{B}_Y$ be the semi-universal deformation space of $X$ and $Y$ respectively. Since $H^2(\mathscr{O}_Y) = 0$, every Enriques manifold is projective and hence its universal cover $X$ is projective as well (see \cite{OS}, Corollary $2.7$). However since non-projective hyperk\"ahler manifolds are dense in the Kuranishi space of $X$, we immediately see that not all deformations of $X$ can be realised as the universal cover of a deformation of $Y$. The elements of $\mathscr{B}_X$ that can be realized as the universal cover of a deformation of $Y$ form a proper closed sublocus $L$ of dimension $ \textrm{dim}(\mathscr{B}_Y) = h^1(T_Y) = h^1(T_X)_0 = h^{1,1}(X)_1$. (Note that $Y = X/G$ where $G$ is cyclic group of $d$-th roots of unity. Hence $G$ acts on both $H^1(T_X)$ and $H^{1,1}(X)$. We denote by $H^1(T_X)_0$ and $H^{1,1}(X)_1$ the summands corresponding to characters identity and a primitive $d$-th root of unity $\zeta$ respectively in the weight decompositions of $H^1(T_X)$ and $H^{1,1}(X)$ indexed by the characters of $G$. See \cite{OS1}, Theorems $1.1, 1.2$ and $2.4$ for details.) \QEDA \par
\end{remark}

We now consider the deformations of the universal covering map of the Hilbert scheme of points on an Enriques surface. Oguiso and Schr\"oer proved in \cite{OS} that the manifold we get is Enriques in our sense since the universal cover is a Calabi-Yau manifold. \par

\begin{theorem}\phantomsection\label{defspcy}
 (\cite{Ha}, Theorem 3.2) Let for $n \geq 2$, $Y=\textrm{Hilb}^n(S)$ be the Hilbert scheme of $n$ points on an Enriques surface. Let $\pi:X \to Y$ be the universal double cover of $Y$. Let $\mathscr{B}_X$ and $\mathscr{B}_{Y}$ be the algebraic (since $H^2(\mathscr{O}_X) = H^2(\mathscr{O}_Y) = 0)$) formally semi-universal deformation spaces of $X$ and $Y$ respectively. Then a general element of $\mathscr{B}_X$ can be realized as the universal double cover of a non-simply connected smooth Enriques variety $Y'$ with $c_1(Y') = 0$, $\chi(\mathscr{O}_{Y'}) = 1$ and $\pi_1(Y')$ cyclic of order $2$.
\end{theorem}   
\noindent\textit{Proof.} We first prove that $H^1(T_Y \otimes K_Y) = h^{1,2n-1}(Y) = 0$. Recall that the Hodge-Poincar\'e polynomial of a smooth projective variety $X$ is defined as $h(X,x,y) = \displaystyle \sum_{p,q =0}^{\dim(X)}h^{p,q}(X)x^py^q$. Using Theorem $2.3.14$ of \cite{Got}, we have that if $S$ denote an Enriques surface with Hodge diamond 
\[
\begin{array}{ccccc}
   & & 1 & &   \\
   & 0 & & 0 & \\
   0 & & 10 & & 0 \\
   & 0 & & 0 & \\
   & & 1 & &  \\
\end{array}
\]
then we have the following generating polynomial for the Hodge numbers of $S^{[n]}$: 
\begin{equation*}
    \displaystyle \sum_{n=1}^{\infty}h(S^{[n]},x,y)t^n = \displaystyle \prod_{k=1}^{\infty}(1-x^{k-1}y^{k-1}t^k)^{-1}(1-x^ky^kt^k)^{-10}(1-x^{k+1}y^{k+1}t^k)^{-1}.
\end{equation*}
Notice that the indices of $x$ and $y$ in each term of the infinite product are the same. Hence we conclude that for each $n$, the resulting Hodge diamond has the property that $h^{p,q}(S^{[n]}) = 0$ if $p \neq q$. Hence $h^{1,2n-1}(S^{[n]}) = 0$ if $n \geq 2$. Since $K_Y$ is the trace zero module of the morphism $\pi$, we have by \cite{Weh}, Corollary 1.11 that the deformation 
 \[
\begin{tikzcd}
    X \arrow{r} \arrow{d} & \mathscr{X} \arrow[d] \\
    0 \arrow{r} & \mathscr{B}_X 
\end{tikzcd} \]
can be lifted to a deformation $\Pi$ \[
\begin{tikzcd}
    X \arrow{r} \arrow[d , "\pi"] & \mathscr{X} \arrow[d , "\Pi"] \\
   Y \arrow{r} \arrow{d} & \mathscr{Y} \arrow{d}  \\
    0 \arrow{r} & \mathscr{B}_X 
\end{tikzcd} \]
of the morphism $\pi$. A general fibre $\mathscr{Y}_t$ of $\mathscr{Y}$ is smooth and non-simply connected with $c_1(\mathscr{Y}_t) = 0$, $\chi(\mathscr{O}_{\mathscr{Y}_t}) = 1$ and $\pi_1(\mathscr{Y}_t)$ cyclic of order $2$. \par 

We can assume that after possibly shrinking $\mathscr{B}_X$ the morphism $\Pi_t$ is finite surjective and is of degree $2$ for all $t \neq 0$. Since $\Pi_t$ is a double cover it is given by a line bundle $B_t$ and a section $s \in H^0(B_t^{\otimes 2})$. But considering the fact that $\mathscr{O}_{\mathscr{X}_t} = K_{\mathscr{X}_t} = \Pi_t^*(K_{\mathscr{Y}_t} \otimes B_t^{-1}) $, we have that $K_{\mathscr{Y}_t} \otimes B_t^{-1}$ is torsion and hence is numerically trivial. Therefore $B_t$ is numerically trivial. Considering that $B_t^{\otimes 2}$ has a section, we conclude that $B_t^{\otimes 2} = \mathscr{O}_X$. Hence $\Pi_t$ is \'etale and by \cite{OS}, Proposition $2.8$ and Lemma $3.2$, $\Pi_t$ is the universal covering morphism and $B_t = K_{\mathscr{X}_t}$. This concludes the proof.\QEDB\par

\vspace{5pt}

We now intend to embed the Enriques manifold inside a projective space and consider the deformations of the composed morphism. This brings up the context of reduced structures, also called ropes on Enriques manifolds. 

\subsection{Ropes on reduced connected schemes.} \phantomsection\label{prelimropes}Ropes are multiple structures on a scheme. Our main objective is to show that ropes on Enriques manifolds can be realized as a flat limit of hyperk\"ahler or Calabi-Yau varieties which we do in section ~\ref{smoothing}. We first recall the definition of ropes.
\begin{definition}\phantomsection\label{defropes}
Let $Y$ be a reduced connected scheme and let $\mathscr{F}$ be a vector bundle of rank $m-1$ on $Y$. A {\it rope of multiplicity $m$ on Y with conormal bundle $\mathscr{F}$} is a scheme $\widetilde{Y}$ with $\widetilde{Y}_{\textrm{red}}=Y$ such that
\begin{itemize}
    \item[(1)] $\mathscr{I}_{Y/\widetilde{Y}}^2=0$,
    \item[(2)] $\mathscr{I}_{Y/\widetilde{Y}}=\mathscr{F}$ as $\mathscr{O}_Y$ modules.
\end{itemize}
If $\mathscr{F}$ is a line bundle then $\widetilde{Y}$ is called a {\it ribbon} on $Y$.
\end{definition}
The following properties of ropes were proven in \cite{BE} and \cite{Gon}.
\begin{theorem}\phantomsection\label{charropes}
Let $Y$ be a reduced connected scheme and $\mathscr{F}$ be a vector bundle of rank $m-1$ on Y.
\begin{itemize}
    \item[(1)] A rope $\widetilde{Y}$ with conormal bundle $\mathscr{F}$ is defined by an element $[e_{\widetilde{Y}}]\in \textrm{Ext}^1(\Omega_Y,\mathscr{F})$. The rope is called split if $[e_{\widetilde{Y}}]=0$.
\end{itemize}
Assume further that $Y$ is a smooth variety and $i:Y\hookrightarrow Z$ is a closed immersion into another smooth variety $Z$.
\begin{itemize}
    \item[(2)] There is an one-to-one correspondence between pairs $(\widetilde{Y},\widetilde{i})$ where $\widetilde{Y}$ is a rope with conormal bundle $\mathscr{F}$ and $\widetilde{i}:\widetilde{Y}\to Z$ is a morphism extending $i:Y\hookrightarrow Z$ and elements $\tau\in \textrm{Hom}(\mathscr{N}_{Y/Z}^*,\mathscr{F})$.
    \item[(3)] If $\tau\in \textrm{Hom}(\mathscr{N}_{Y/Z}^*,\mathscr{F})$ corresponds to $(\widetilde{Y},\widetilde{i})$, then $\widetilde{i}$ is an embedding if and only if $\tau$ is surjective.
    \item[(4)] If $\tau\in \textrm{Hom}(\mathscr{N}_{Y/Z}^*,\mathscr{F})$ corresponds to $(\widetilde{Y},\widetilde{i})$, then $\tau$ is mapped by the connecting homomorphism onto $[e_{\widetilde{Y}}]$.
\end{itemize}
\end{theorem}
We will encounter $d$-ropes on Enriques manifolds of index $d$. Of these, the $2$-ropes -- also called ribbons -- have pleasant geometric properties and in our context, they are analogues of $K3$-carpets. We recall some of their properties. It is known that a ribbon on a smooth, irreducible variety $Y$ is locally Gorenstein and consequently it will have a dualizing sheaf $K_{\widetilde{Y}}$ which is a line bundle (see \cite{GGP}, Remark 1.3). In \cite{GGP}, the authors defined $K3$ carpets on a regular surface to be a ribbon with trivial canonical bundle. We generalize their definition in the following way to define $K$- trivial ribbons. The name is pretty self-explanatory.
\begin{definition}\phantomsection\label{defcyr}
Let $Y$ be a smooth, irreducible, projective variety, $\widetilde{Y}$ be a ribbon on $Y$ with conormal bundle $\mathscr{F}$. Then {\it $\widetilde{Y}$ is a $K$- trivial ribbon} if $K_{\widetilde{Y}}=\mathscr{O}_{\widetilde{Y}}$. 
\end{definition}
A ribbon $\widetilde{Y}$ is regular if $h^1(\mathscr{O}_{\widetilde{Y}})=0$. We include a characterization of a $K$- trivial ribbon from its conormal bundle whose proof is similar to that of \cite{GGP}, Proposition 1.5.
\begin{lemma}\phantomsection\label{charcyr}
Let Y be a smooth, irreducible, projective variety of dimension n. Let $h^{n-1}(\mathscr{O}_Y)=0$. Let $\widetilde{Y}$ be a ribbon with conormal bundle $\mathscr{F}$. Then $\widetilde{Y}$ is $K$- trivial if and only if $\mathscr{F}\cong K_Y$. Moreover, if $h^{n-2}(\mathscr{O}_Y)=0$ then $H^1(\mathscr{O}_{\widetilde{Y}})\cong H^1(\mathscr{O}_Y)$, if $Y$ is regular then $\widetilde{Y}$ is also regular.
\end{lemma}
\noindent\textit{Proof.} Notice that, we have the following short exact sequence;
\begin{equation*}
    0\to\mathscr{F}\to \mathscr{O}_{\widetilde{Y}}\to\mathscr{O}_Y\to 0.
\end{equation*}
It follows from Lemma 1.4, (3) \cite{GGP} that if we apply  $\mathscr{H}om(\textrm{--},K_{\widetilde{Y}})$ on the above exact sequence, we get the following exact sequence;
\begin{equation*}
    0\to K_Y\to K_{\widetilde{Y}}\to\mathscr{F}^{-1}\otimes K_Y\to 0.
\end{equation*}
\indent First, assume that $\mathscr{F}=K_Y$. It follows from the previous exact sequence that we have a surjection $H^0(K_{\widetilde{Y}})\to H^0(\mathscr{O}_Y)$ since $h^1(K_Y)=0$ by assumption. Thus $K_{\widetilde{Y}}=\mathscr{O}_{\widetilde{Y}}$. \par Now, assume $K_{\widetilde{Y}}=\mathscr{O}_{\widetilde{Y}}$. We tensor the above exact sequence by $\mathscr{O}_Y$ to get a surjection $\mathscr{O}_Y\to\mathscr{F}^{-1}\otimes K_Y$ and hence $\mathscr{F}=K_Y$.\par  The remaining assertions are clear from the first exact sequence.\QEDB\par 

\vspace{5pt}

We will show that $d$-ropes with conormal bundle $\mathscr{E}$, appear naturally on Enriques manifolds of index $d$ as flat limits of hyperk\"ahler or Calabi-Yau manifolds. Our proofs are aided by the following main theorem that connects deformation of finite morphisms to ropes on their images. \par
Suppose $X$ and $Y$ are smooth, irreducible projective varieties. Let $i: Y \hookrightarrow \mathbb{P}^N$ be an embedding. Let $\mathscr{I}$ be the ideal sheaf of $i(Y)$ inside $\mathbb{P}^N$. Let $\pi: X \to Y$ denote a finite morphism of degree $d \geq 2$ with trace zero module given by $\mathscr{F}$ where $\mathscr{F}$ is a vector bundle on $Y$ of rank $d-1$. Let $\varphi = \pi \circ i$. We know that $H^0(\mathscr{N}_{\varphi})$ parametrizes the first order deformations of $\varphi$. Let $\Psi_2: H^0(\mathscr{N}_{\varphi}) \to \textrm{Hom}(\mathscr{I}/\mathscr{I}^2,\mathscr{F})$ be as in Proposition $1.2$, \cite{GGP1}.
Assuming this setting, the following result has been proven by Gallego, Gonz\'alez and Purnaprajna.
\begin{theorem}\phantomsection\label{defs and ropes}
Let $T$ be a smooth irreducible algebraic curve and let $0$ be a closed point of $T$. Let $\Phi: \mathscr{X} \to \mathbb{P}_T^N$ be a flat family of morphisms over $T$ (i.e, $\Phi$ is a $T-$ morphism such that $\mathscr{X} \to T$ is flat, proper and surjective) such that 
\begin{itemize}
    \item[(a)] $\mathscr{X}$ is reduced and irreducible
    \item[(b)] $\mathscr{X}_t$ is smooth, irreducible and projective for all $t \in T$
    \item[(c)] $\mathscr{X}_0 = X$ and $\Phi_0 = \varphi$.
\end{itemize}
Let $\Delta$ be the first infinitesimal neighbourhood of $0 \in T$. Let $\widetilde{X} = \mathscr{X}_{\Delta} := \mathscr{X} \times_T \Delta$ and $\widetilde{\varphi} = \Phi_{\Delta}: \widetilde{X} \to \mathbb{P}_{\Delta}^N$ be the first order deformations of $X$ and $\varphi$ respectively induced by the morphism $\Delta \to T$. Let $\nu$ be the element in $H^0(\mathscr{N}_{\varphi})$ corresponding to $\widetilde{\varphi}$.
\begin{itemize}
    \item[(1)] If $\Psi_2(\nu) \in \textrm{Hom}(\mathscr{I}/\mathscr{I}^2,\mathscr{F})$ is a homomorphism of rank $> d/2 - 1$, then after shrinking $T$ if necessary, $\Phi_t$ is one-to-one for all $t \in T, t \neq 0$.
   
    \item[(2)] If $\Psi_2(\nu) \in \textrm{Hom}(\mathscr{I}/\mathscr{I}^2,\mathscr{F})$ is a surjective homomorphism, then after shrinking $T$ if necessary, $\Phi_t$ is a closed immersion for all $t \in T, t \neq 0$. In this case $(\textrm{Im}\, \Phi)_0 = (\textrm{Im} \,\widetilde{\varphi})_0$ is an embedded rope on $Y$ with conormal bundle $\mathscr{F}$ .
\end{itemize}

\end{theorem}

\noindent\textit{Proof.} For part $(1)$ see \cite{GGP2}, Proposition $1.4$. For part $(2)$, see \cite{GGP1}, Proposition $1.4$. The last assertion of part $(2)$ follows from equation $1.9$ in Proposition $1.4$, \cite{GGP2}. \QEDB 

\section{Ropes supported on Enriques manifolds}\phantomsection\label{RE}
Recall that for a given morphism $\varphi:X\to \mathbb{P}^N$ from a variety $X$ to projective space, $H^0(\mathscr{N}_{\varphi})$ parametrizes the first order infinitesimal deformation of $\varphi$ and $H^1(\mathscr{N}_{\varphi})$ is its obstruction space. In \cite{GGP}, the authors proved that $\varphi$ is unobstructed for a $K3$ surface $X$ if $\varphi$ is finite onto its image, by showing $h^1(\mathscr{N}_{\varphi})=0$. However, in our case $h^1(\mathscr{N}_{\varphi})\neq0$ in general but $\varphi$ will still be unobstructed as the next two lemmas show.
\begin{lemma}\phantomsection\label{unobofphi}
Let $X$ be a hyperk\"ahler variety and let $X\xrightarrow{\varphi}\mathbb{P}^N$ be a morphism which is finite onto its image. Then $\varphi$ is unobstructed and $H^1(\mathscr{N}_{\varphi})\cong H^2(T_X)$.
\end{lemma}
\noindent\textit{Proof.} We start with the Atiyah extension of $L=\varphi^*(\mathscr{O}_{\mathbb{P}^N}(1))$ (see \cite{Ser})
\begin{align}\phantomsection\label{AE}
    0 \rightarrow \mathscr{O}_X\rightarrow \mathscr{E}_L\rightarrow T_X\rightarrow 0.
\end{align}
Since $H^1(\mathscr{O}_X)=0$, the long exact sequence associated with this short exact sequence gives
\begin{align}\phantomsection\label{LEA}
    0 \rightarrow H^1(\mathscr{E}_L)\rightarrow H^1(T_X)\rightarrow H^2(\mathscr{O}_X).
\end{align}
Recall that the map $H^1(T_X)\to H^2(\mathscr{O}_X)$ is induced by $c_1(L)\in H^1(\Omega_X)$ and the natural pairing $H^1(T_X)\otimes H^1(\Omega_X)\to H^2(\mathscr{O}_X)$ is perfect (see \cite{Hu}, 1.8). Therefore the rightmost map in ~\ref{LEA} is surjective and consequently $h^1(\mathscr{E}_L)=h^1(T_X)-1$.\par 
\indent Next, consider the following exact sequence 
\begin{align}\phantomsection\label{TTN}
    0\to T_X\to \varphi^*T_{\mathbb{P}^N}\to \mathscr{N}_{\varphi}\to 0.
\end{align}
The long exact sequence associated to this short exact sequence gives 
\begin{align}\phantomsection\label{LETTN}
    H^0(\mathscr{N}_{\varphi})\xrightarrow{\nu} H^1(T_X)\xrightarrow{} H^1(\varphi^*T_{\mathbb{P}^N})\xrightarrow{} H^1(\mathscr{N}_{\varphi})\xrightarrow{\psi} H^2(T_X).
\end{align}
Note that $\nu$ factors through $H^1(\mathscr{E}_L)$. Thus $\nu$ is not surjective since $h^1(\mathscr{E}_L)<h^1(T_X)$. \par
Using the pull back of the Euler sequence on $\mathbb{P}^N$, it can be seen easily that $h^1(\varphi^*T_{\mathbb{P}^N})=1$, $h^2(\varphi^*T_{\mathbb{P}^N})=0$. Consequently, $\psi$ is bijective and the second assertion follows.\par 
Recall that $H^1(\mathscr{N}_{\varphi})$ is an obstruction space of ${\bf Def}_{\varphi}$ (deformation of $\varphi$ with fixed target) and $H^2(T_X)$ is an obstruction space of ${\bf Def}_X$ (deformation of $X$). The fact that $\varphi$ is unobstructed follows from the fact that $X$ is unobstructed (see \cite{Ti}, \cite{To}).\QEDB

\vspace{5pt}

In the situation of the lemma above, if we had an even dimensional Calabi-Yau manifold $X$ instead of the hyperk\"ahler manifold, the long exact sequence of the restricted Euler sequence would have showed $h^i(\varphi^*T_{\mathbb{P}^N})=0$ for $i=1,2$. Consequently the long exact sequence ~\eqref{LETTN} would again show bijection between $H^1(\mathscr{N}_{\varphi})$ and $H^2(T_X)$ and the same argument would yield the unobstructedness of $\varphi$. We formulate the discussion in the following lemma.
\begin{lemma}\phantomsection\label{unobphicy}
Let $X$ be a Calabi-Yau manifold and let $X\xrightarrow{\varphi}\mathbb{P}^N$ be a morphism which is finite onto its image. Then $\varphi$ is unobstructed and $H^1(\mathscr{N}_{\varphi})\cong H^2(T_X)$.
\end{lemma}
In \cite{GGP} the authors proved that for an Enriques surfaces $Y$ and a very ample line bundle $L$ on $Y$, the line bundle $L\otimes K_Y$ is globally generated and the pull back of $L$ to the universal $K3$ cover is very ample using Reider's theorem which does not carry over to higher dimension. However, we show that the pull back of a very ample line bundle $L$ on an Enriques manifold remains very ample on its universal cover provided $L\otimes K_Y^{\otimes i}$ is globally generated for all $1\leq i\leq d-1$. This result might be of independent interest even though will will never use it in the proof of our subsequent results.

\begin{lemma}\phantomsection\label{va}
Let $Y$ be an Enriques manifold of index $d$ and let $\pi:X\to Y$ be its universal cover. Suppose $L$ is a very ample line bundle on $Y$ for which $L\otimes K_Y^{\otimes i}$ is globally generated for all $1\leq i\leq d-1$. Then $\pi^*L$ is very ample.
\end{lemma}
\noindent\textit{Proof.} Since $L$ is very ample, certainly $\pi^*L$ separates two points that belong to different fibers. Moreover, $\pi^*L$ separates tangent vectors since \'etale morphism induces isomorphism of tangent spaces. Therefore, $\pi^*L$ is very ample $\iff$ $\pi^*L$ separates two points that belong to the same fiber.\par 
\indent Let $\zeta$ be the length $d$ subscheme supported on $p_i\in X$ for $1\leq i\leq d$ that maps to a point $q\in Y$. We have the following short exact sequence on $Y$
\begin{align}
    0\to \mathscr{I}_q\to\mathscr{O}_Y\to\mathscr{O}_q\to 0.
\end{align}
Pulling it back on $X$ and tensoring by $\pi^*L$ we get
\begin{align}
    0\to \pi^*\left(L\otimes\mathscr{I}_q\right)\to\pi^*L\to \pi^*L\otimes\mathscr{O}_{\zeta}\to 0.
\end{align}
Note that $h^1(L\otimes K_Y^{\otimes i}\otimes\mathscr{I}_q)=0$ for all $i$ by our assumption. Consequently, $h^1(\pi^*(L\otimes\mathscr{I}_q))=0$ which can be seen by pushing forward $\pi^*(L\otimes\mathscr{I}_q)$ (projection formula applies since $\pi$ is finite, see \cite{A}, Lemma 5.7). Thus, $H^0(\pi^*L)$ separates points of a fiber.\QEDB 

\vspace{5pt}

The following proposition enables us to show the existence of surjective homomorphisms in the group $\textrm{Hom}(\mathscr{I}/\mathscr{I}^2,\mathscr{E}) = H^0(\mathscr{N}_{Y/\mathbb{P}^N}\otimes \mathscr{E})$ 
for an embedding of an Enriques manifold of index $d$ in a large enough projective space. \par

Before stating the result, let us introduce the invariant which we call the sectional genus $g$ and which by definition is the dimension of the vector space of global sections of the canonical bundle of a smooth hyperplane section of the embedding. Since $H^1(K_Y) = 0$, one can see that $g = h^0(K_Y(1))$. Also since $K_Y$ is numerically trivial and $H^j(K_Y^{\otimes i}(1)) = 0$ for $1 \leq i \leq d-1$ and $j > 0$, we have that $g = h^0(K_Y^{\otimes i}(1))$ for $1 \leq i \leq d-1$.

\begin{proposition}\phantomsection\label{countinghk}
Let $Y$ be an Enriques manifold of and index $d$ and dimension $2n$. Let $X$ be the hyperk\"ahler universal cover of $Y$. Let $i:Y\hookrightarrow\mathbb{P}^N$ be an embedding of $Y$. Then,
\begin{itemize}
    \item[(1)] The non-split abstract multiple structures on $Y$ with conormal bundle $\mathscr{E}$ are parametrized by the projective space of lines of the vector space $H^1(T_Y \otimes \mathscr{E})$ of dimension $h^1(T_X)-h^1(T_Y)$. Among them, the projective multiple structures form a countable union of hyperplanes.  
    \item[(2)] The dimension of the vector space $H^0(\mathscr{N}_{Y/\mathbb{P}^N}\otimes \mathscr{E})$ that parametrizes pairs $(\widetilde{Y}, \widetilde{i})$ where $\widetilde{Y}$ is a rope on $Y$ of multiplicity $d$ with conormal bundle $\mathscr{E}$ and $\widetilde{i}$ is a morphism extending $i$ is given by 
\begin{align*}
    h^0(\mathscr{N}_{Y/\mathbb{P}^N}\otimes\mathscr{E})=(d-1)g(N+1)+h^1(T_X)-h^1(T_Y)-1
\end{align*}
where $g$ is the sectional genus of the embedding.
    \item[(3)] Suppose $N \geq 4n+d-1$ and $\mathscr{O}_Y(1)\otimes K_Y^{\otimes i}$ is globally generated for all $1\leq i\leq d-1$. 
Then the $d$-rope structures with conormal bundle $\mathscr{E}$ embedded in $\mathbb{P}^N$ and supported on $i(Y)$ are parametrized by a non-empty open set in the projective space of lines in the space $H^0(\mathscr{N}_{Y/\mathbb{P}^N}\otimes \mathscr{E})$. 
\end{itemize}
\end{proposition} 

\noindent\textit{Proof.}  ($1$) By Theorem ~\ref{charropes}, we have that the abstract ropes with conormal bundle $\mathscr{E}$ are characterized by $\textrm{Ext}^1(\Omega_Y, \mathscr{E}) = H^1(T_Y \otimes \mathscr{E})$ since $Y$ is smooth. Now we follow the proof of \cite{GGP}, Theorem $2.5, (ii)$. Since the ideal of $Y$ inside $\widetilde{Y}$ is square-zero, given a multiple structure $\widetilde{Y}$ we have an exact sequence
\begin{equation*}
    0 \to \mathscr{E} \to \mathscr{O}_{\widetilde{Y}}^* \to \mathscr{O}_Y^* \to 1
\end{equation*}
which induces the following exact sequence
\begin{equation}\label{pic}
   \textrm{Pic}(\widetilde{Y}) \to \textrm{Pic}(Y) \to H^2(\mathscr{E}) = \mathbb{C}. 
\end{equation}
A line bundle $\widetilde{A} \in \textrm{Pic}(\widetilde{Y})$ is ample if and only if its restriction $A \in \textrm{Pic}(Y)$ is ample. Consider the universal cover $X$. We have a non-degenerate cup-product pairing
\begin{equation*}
H^1(\Omega_X) \otimes H^1(T_X) \to H^2(\mathscr{O}_X).
\end{equation*}
Pushing forward we have a map
\begin{equation*}
(H^1(\Omega_Y) \oplus H^1(\Omega_Y \otimes \mathscr{E})) \otimes (H^1(T_Y) \oplus H^1(T_Y \otimes \mathscr{E})) \to H^2(\mathscr{O}_Y) \oplus H^2(\mathscr{E}).  
\end{equation*}
Since $H^2(\mathscr{O}_Y) = 0$ and the pairing on $X$ is non-degenerate, we have the pairing
\begin{equation*}
    H^1(\Omega_Y) \otimes H^1(T_Y \otimes \mathscr{E}) \to H^2(\mathscr{E}) = \mathbb{C}
\end{equation*}
is non-degenerate. Since the ample cone of $Y$ is contained inside $H^1(\Omega_Y)$, we have that for each ample line bundle $A$ in $\textrm{Pic}(Y)$, there exists a hyperplane in $\mathbb{P}(H^1(T_Y \otimes \mathscr{E}))$ such that $A$ is in the kernel of  the map $\textrm{Pic}(Y) \to H^2(\mathscr{E}) = \mathbb{C}$ in ~\eqref{pic}. Hence the projective multiple structures form a countable union of hyperplanes. \par 

(2) Let us first compute  $ h^0(\mathscr{N}_{Y/\mathbb{P}^N} \otimes \mathscr{E}) = \sum\limits_{i=1}^{d-1} h^0(\mathscr{N}_{Y/\mathbb{P}^N} \otimes K_Y^{\otimes i})$. 
We start with the following short exact sequence
\begin{align}\phantomsection\label{ttn}
    0\to T_Y\to T_{\mathbb{P}^N}\vert_Y\to\mathscr{N}_{Y/\mathbb{P}^N}\to 0.
\end{align}
Let $i$ be an integer in between 1 and $d-1$. Note that $h^0(T_Y\otimes K_Y^{\otimes i})=0$. We tensor the above exact sequence by $K_Y^{\otimes i}$ and take the long exact sequence of cohomology to get the following. 
\begin{align}\phantomsection\label{le1c}
\begin{split}
    0\to H^0(T_{\mathbb{P}^N}\vert_Y\otimes K_Y^{\otimes i})\to H^0(\mathscr{N}_{Y/\mathbb{P}^N}\otimes K_Y^{\otimes i}) \to H^1(T_Y\otimes K_Y^{\otimes i})\xrightarrow{\alpha} H^1(T_{\mathbb{P}^N}\vert_Y\otimes K_Y^{\otimes i})\\ \to H^1(\mathscr{N}_{Y/\mathbb{P}^N}\otimes K_Y^{\otimes i})\to H^2(T_Y\otimes K_Y^{\otimes i})\to H^2(T_{\mathbb{P}^N}\vert_Y\otimes K_Y^{\otimes i})\to ...
\end{split}    
\end{align}
Restrict of the Euler sequence of $\mathbb{P}^N$ on $Y$ and tensor it by $K_Y^{\otimes i}$ to obtain the following exact sequence
\begin{align}
    0\to K_Y^{\otimes i}\to K_Y^{\otimes i}(1)^{\oplus N+1}\to T_{\mathbb{P}^N}\vert_Y\otimes K_Y^{\otimes i}\to 0.
\end{align}
Using the cohomology sequence associated to the short exact sequence above one can check that
\begin{align}
    h^0(T_{\mathbb{P}^N}\vert_Y \otimes K_Y^{\otimes i}) = g(N+1),\quad 
    h^1(T_{\mathbb{P}^N}\vert_Y \otimes K_Y^{\otimes i}) \leq 1,\quad
    h^2(T_{\mathbb{P}^N}\vert_Y \otimes K_Y^{\otimes i}) = 0.
\end{align}
 Also, since $h^2(\mathscr{O}_X) = 1$ we have that there exists exactly one $i$ such that $h^1(T_{\mathbb{P}^N}\vert_Y \otimes K_Y^{\otimes i}) = 1$. \par 
Let $\pi:X\to Y$ be the universal cover and $\varphi=i\circ\pi$. Since $\pi$ is \'etale, $\mathscr{N}_{\varphi}\cong \pi^*\mathscr{N}_{Y/\mathbb{P}^N}$. Consequently we have the following,
\begin{align}
    h^1(\mathscr{N}_{\varphi})= \sum_{i=0}^{d-1}h^1(\mathscr{N}_{Y/\mathbb{P}^N}\otimes K_Y^{\otimes i}).
\end{align}
By Lemma ~\ref{unobofphi} we also have that $h^1(\mathscr{N}_{\varphi}) = h^2(T_X)$ and consequently we get,
\begin{align}
    h^1(\mathscr{N}_{\varphi}) = h^2(T_X) = \sum_{i=0}^{d-1}h^2(T_Y \otimes K_Y^{\otimes i}).
\end{align}
Comparing the above two equations one can see at once that
\begin{align}
    \sum\limits_{i=1}^{d-1}h^1(\mathscr{N}_{Y/\mathbb{P}^N} \otimes K_Y^{\otimes i}) - \sum\limits_{i=1}^{d-1}h^2(T_Y \otimes K_Y^{\otimes i}) = h^2(T_Y)-h^1(\mathscr{N}_{Y/\mathbb{P}^N}) = 0
\end{align}
where the last equality is obtained by the restriction of the Euler sequence of $\mathbb{P}^N$ to $Y$. Also notice that since $T_X = \pi^*(T_Y)$ we have that $h^1(T_X)=h^1(T_Y)+\sum\limits_{i=1}^{d-1}h^1(T_Y\otimes K_Y^{\otimes i})$.\par
\indent Using the long exact sequence ~\eqref{le1c} we conclude that 
\begin{align*}
    \sum_{i=1}^{d-1}h^0(\mathscr{N}_{Y/\mathbb{P}^N}\otimes K_Y^{\otimes i})
= (d-1)g(N+1) + \sum_{i=1}^{d-1}h^1(T_Y \otimes K_Y^{\otimes i}) - \sum_{i=1}^{d-1}h^2(T_Y \otimes K_Y^{\otimes i})\\ + \sum_{i=1}^{d-1}h^1(\mathscr{N}_{Y/\mathbb{P}^N} \otimes K_Y^{\otimes i}) - \sum_{i=1}^{d-1}h^1(T_{\mathbb{P}^N}\vert_Y \otimes K_Y^{\otimes i})-1\\ = (d-1)g(N+1)+h^1(T_X)-h^1(T_Y)-1.
\end{align*}

(3) Recall that embedded multiple structures with conormal bundle  $\mathscr{E}$ are in one-to-one correspondence with surjective homomorphisms in $\textrm{Hom}(\mathscr{I}/\mathscr{I}^2, \mathscr{E})$ upto non-zero scalar multiple. \par 
Now, for a morphism $\nu:\mathscr{F}\to \mathscr{G}$ between vector bundles of rank $f$ and $g$ respectively, on an irreducible complex projective  variety $W$, for any positive integer $k\leq \textrm{min}\{f,g\}$, the {\it $k$-th degeneracy locus $D_k(\nu)$} is defined as the subscheme cut out by the minors of order $k+1$ of the matrix locally representing $\nu$. A result of \cite{B} asserts that if $\mathscr{F}^*\otimes\mathscr{G}$ is globally generated, then for a general $\nu$, either $D_k(\nu)$ is empty, or has pure codimension $(f-k)(g-k)$.\par

Note that we have a surjection $T_{\mathbb{P}^N}|_Y \otimes \mathscr{E} \longrightarrow \mathscr{N}_{Y/\mathbb{P}^N}\otimes\mathscr{E}$ and by the Euler sequence on $\mathbb{P}^N$ twisted by $\mathscr{E}$, we conclude that $\mathscr{N}_{Y/\mathbb{P}^N}\otimes\mathscr{E}$ globally generated since $\mathscr{O}_Y(1)\otimes \mathscr{E} $ is globally generated by our assumptions. Now we apply the result explained in the above paragraph with $\mathscr{F}=\mathscr{I}/\mathscr{I}^2$ and $\mathscr{G}=\mathscr{E}$ (recall that $\mathscr{N}_{Y/\mathbb{P}^N}=(\mathscr{I}/\mathscr{I}^2)^*$). Thus, for a general homomorphism in $\textrm{Hom}(\mathscr{I}/\mathscr{I}^2, \mathscr{E})$, $D_{d-2}$ is either empty, or of pure codimension $N-\dim(Y)-(d-2)$. However, the latter is impossible since by our assumption, $N-\dim(Y)-(d-2)>\dim(Y)$.\QEDB\par

\vspace{5pt}
 
We can perform the analogous computations for an Enriques manifold whose universal cover is a Calabi-Yau manifold. Note that in this case $d$ is always equal to $2$ and $h^1(T_X)-h^1(T_Y)=h^{1,2n-1}(Y)$ where $\dim(Y)=2n$. Also note that in this case all multiple structures with conormal bundle $\mathscr{E}$ are projective since $H^2(K_Y) = 0$. The result we obtained is given below whose proof we omit.
\begin{proposition}\phantomsection\label{countingcy}
Let $Y$ be an Enriques manifold of dimension 2n whose universal cover is a Calabi-Yau manifold. Let $i:Y\hookrightarrow\mathbb{P}^N$ be an embedding of $Y$. 
\begin{itemize}
    \item[(1)] The non-split abstract multiple structures with conormal bundle $\mathscr{E}$ supported on $Y$ are parametrized by the projective space of lines of the vector space $H^{2n-1}(\Omega_Y)^*$. All such multiple structures are projective.
    \item[(2)] The dimension of the vector space $H^0(\mathscr{N}_{Y/\mathbb{P}^N}\otimes \mathscr{E})$ that parametrizes pairs $(\widetilde{Y}, \widetilde{i})$ where $\widetilde{Y}$ is a rope on $Y$ of multiplicity $2$ with conormal bundle $\mathscr{E}$ and $\widetilde{i}$ is a morphism extending $i$ is given by 
\begin{align*}
    h^0(\mathscr{N}_{Y/\mathbb{P}^N}\otimes \mathscr{E})=g(N+1)+h^{1,2n-1}(Y)
\end{align*}
where $g$ is the sectional genus of the embedding.
\item[(3)] Assume that $N \geq 4n+1$ and that $\mathscr{O}_Y(1)\otimes K_Y$ is globally generated. Then the double structures with conormal bundle $K_Y$ embedded in $\mathbb{P}^N$ and supported on $i(Y)$ are parametrized by a non-empty open set in the projective space of lines in $H^0(\mathscr{N}_{Y/\mathbb{P}^N}\otimes K_Y)$. 
\end{itemize}
\end{proposition}
Now we discuss the condition of global generation of $\mathscr{O}_Y(1) \otimes K_Y^{\otimes i}$ for $1 \leq i \leq d-1$. In particular we discuss how we can get a very ample line bundle $L$ such that $L \otimes K_Y^{\otimes i}$ is base point free for $1 \leq i \leq d-1$.
\begin{remark}
 Let $Y$ be an Enriques manifold of dimension $2n$ and $L$ be an ample line bundle on $Y$.
 \begin{itemize}
     \item[(1)] Notice that if we start with an ample line bundle $L$, there exists $M$ depending on $L$ such that $L^{ \otimes m}$ is very ample and $L^{\otimes m} \otimes K_Y^{\otimes i}$ is base point free for $1 \leq i \leq d-1$ and $m \geq M$. 
     \item[(2)] If we start with an ample and base point line bundle $L$, we have by \cite{Pop}, Theorem $1.4$ that $L^{\otimes m} \otimes K_Y^{\otimes i}$ is base point free for all $1\leq i\leq d-1$ if we choose $m \geq d(2n+1)$ where $d$ is the index of $Y$. Also by Castelnuovo-Mumford regularity we have that $L^{\otimes m}$ is very ample if $m \geq 2n+1$. Now since $d \geq 2$, we see that choosing $m \geq d(2n+1)$ is enough to satisfy the conditions.\QEDA
 \end{itemize}
 \end{remark}
\section{Smoothing and its consequence}\phantomsection\label{smoothing}
This section is devoted to prove the existence of an embedded smoothing of projective ropes of multiplicity $d$ with conormal bundle $\mathscr{E} = \bigoplus\limits_{i=1}^{d-1}K_Y^{\otimes i}$ embedded inside a projective space and supported on an Enriques manifold $Y$ of index $d$ provided they exist. Notice that if an embedding $i:Y\hookrightarrow \mathbb{P}^N$ of the Enriques manifold of index $d$ satisfies $N\geq 2\dim(Y)+d-1$ and $\mathscr{O}_Y(1) \otimes K_Y^{\otimes i}$ is globally generated for $1\leq i\leq d-1$, then such a rope exists by Propositions ~\ref{countinghk}, ~\ref{countingcy}. As a consequence of smoothing, we will prove the smoothness of the point corresponding to the rope in its Hilbert scheme when $d=2$. We have included a proof of the fact that the condition $N\geq 2\dim(Y)+1$ is satisfied when the index is two and the universal cover is a hyperk\"ahler six-fold of one of the known deformation types in the Appendix.\par 
First, we will prove the existence of smoothing when the universal cover $X$ is hyperk\"ahler. Our results on smoothing are obtained by using the deformation theoretic results in Theorem ~\ref{defs and ropes}. 

\begin{theorem}\phantomsection\label{smoothinghk}
Suppose $Y$ be an Enriques manifold of index $d$ and dimension $2n$ and $\pi : X \to Y$ be the hyperk\"ahler universal cover of $Y$. Let $i:Y\hookrightarrow\mathbb{P}^N$ be an embedding of $Y$ into a projective space. Let $\varphi = i \circ \pi$.
\begin{itemize}
\item[(i)] Assume either (a) or (b) holds:
\begin{enumerate}
    \item[(a)] $2\leq d\leq 3$; or
    \item[(b)] $d\geq 4$, $\mathscr{O}_Y(1)\otimes K_Y^{\otimes i}$ is globally generated for $1\leq i\leq d-1$, and $$N-2n-\floor{d/2}+1>0,$$
\end{enumerate} 
then there exists a smooth irreducible family $\mathscr{X}$ proper and flat over a smooth pointed \color{black} affine algebraic curve $(T,0)$ \color{black}and a $T$-morphism $\Phi:\mathscr{X} \to \mathbb{P}^N_T $ with the following properties:
\begin{itemize}
    \item[(1)] the fibres $\Phi_t:\mathscr{X}_t \to \mathbb{P}^N$ are finite birational morphisms for $t \in T-\{0\}$,
    \item[(2)] the central fibre $\Phi_0:\mathscr{X}_0 \to \mathbb{P}^N$ is the morphism $\varphi:X \to \mathbb{P}^N$. 
\end{itemize}
\item[(ii)] Let $ \varphi $ be $i \circ \pi $ and $\widetilde{Y}\hookrightarrow\mathbb{P}^N$ be an embedded multiplicity $d$-rope with conormal bundle $\mathscr{E} = \bigoplus\limits_{i=1}^{d-1}K_Y^{\otimes i}$ (which exists under the assumptions of Proposition ~\ref{countinghk}, (3)). Then there exists a smooth irreducible family $\mathscr{X}$ proper and flat over a smooth pointed \color{black} affine algebraic curve $(T,0)$ \color{black}and a $T$-morphism $\Phi:\mathscr{X} \to \mathbb{P}^N_T $ with the following properties:
\begin{itemize}
    \item[(1)] the fibres $\Phi_t:\mathscr{X}_t \to \mathbb{P}^N$ are closed immersions of smooth hyperk\"ahler varieties for $t \in T-\{0\}$,
    \item[(2)] the central fibre $\Phi_0:\mathscr{X}_0 \to \mathbb{P}^N$ is the morphism $\varphi:X \to \mathbb{P}^N$,
    \item[(3)] $\textrm{Im}\,\Phi$ is a flat family of schemes over $T$ embedded inside $\mathbb{P}^N_T$ such that for $t \neq 0$, $(\textrm{Im}\,\Phi)_t$ is a smooth hyperk\"ahler variety and the central fibre $(\textrm{Im}\,\Phi)_0$ is $\widetilde{Y}$.
\end{itemize}
\end{itemize}
\end{theorem} 

\noindent\textit{Proof.} To start with, notice that under the assumptions of part (i) (a), there is a homomorphism of rank $\geq 1$ in $\textrm{Hom}(\mathscr{I}/\mathscr{I}^2,\mathscr{E})$, thanks to Proposition ~\ref{countinghk}. Furthermore, under the assumptions of part (i) (b), a general homomorphism in $\textrm{Hom}(\mathscr{I}/\mathscr{I}^2,\mathscr{E})$ has rank $>\floor{d/2}-1$. Indeed, the result of \cite{B} used in the proof of Proposition ~\ref{countinghk} (3), we obtain that $D_{\floor{d/2}-1}(\nu)$ is a proper closed subvariety of $Y$ by our assumption. Also recall that the given $\widetilde{Y}\hookrightarrow\mathbb{P}^N$ in part (ii) corresponds to a surjective homomorphisms in $\textrm{Hom}(\mathscr{I}/\mathscr{I}^2,\mathscr{E})$ (see Theorem ~\ref{charropes} ($3$)). \par 

\noindent\underline{\textit{Step 1.}} In this step we show the existence of first order deformations  $\widetilde{\varphi}: \widetilde{X} \to \mathbb{P}_{\Delta}^N$ of $\varphi: X \to \mathbb{P}^N$ such that their class $\widetilde{\varphi} \in H^0(\mathscr{N}_{\varphi})$ maps to a rank $\geq 1$ (resp. rank $>\floor{d/2}-1$) homomorphism in part (i) (a) (resp. part (i) (b)) or to the class of $\widetilde{Y}\hookrightarrow\mathbb{P}^N$ in $\textrm{Hom}(\mathscr{I}/\mathscr{I}^2,\mathscr{E})$ in part (ii).\par  

\noindent\textit{Proof of Step 1.} Consider the following short exact sequence 
\begin{align}
    0 \to \mathscr{N}_{\pi} \to \mathscr{N}_{\varphi} \to \mathscr{H}om\left(\pi^*\frac{\mathscr{I}}{\mathscr{I}^2}, \mathscr{O}_X\right) \to 0.
\end{align}
The long exact sequence of cohomology associated to the  sequence above is the following 
\begin{align}
    0 \to H^0(\mathscr{N}_{\pi}) \to H^0(\mathscr{N}_{\varphi})\to \textrm{Hom}\left(\frac{\mathscr{I}}{\mathscr{I}^2}, \mathscr{O}_Y\right)  \oplus \textrm{Hom}\left(\frac{\mathscr{I}}{\mathscr{I}^2}, \mathscr{E}\right)  \to H^1(\mathscr{N}_{\pi}).
\end{align}
Now, $\mathscr{N}_{\pi}$ is supported on the ramification locus
and $\pi$ is an \'etale morphism, consequently we have that $H^1(\mathscr{N}_{\pi}) = 0$. That proves the assertion.\par 

\noindent\underline{\textit{Step 2.}} In this step we show that given a first order deformation $\widetilde{\varphi}: \widetilde{X} \longrightarrow \mathbb{P}_{\Delta}^N$ of $\varphi: X \longrightarrow \mathbb{P}^N$,  there exists a smooth irreducible family $\mathscr{X}$ proper and flat over a 
smooth pointed affine algebraic curve $(T,0)$ and a $T$-morphism $\Phi:\mathscr{X} \to\mathbb{P}^N_T $ such that all squares in the following diagram are Cartesian  \[
\begin{tikzcd}
    X \arrow{r} \arrow[d, "\varphi"] & \widetilde{X} \arrow{r} \arrow[d, "\widetilde{\varphi}"] & \mathscr{X} \arrow{d}\\
   \mathbb{P}^N \arrow{r} \arrow{d} & \mathbb{P}_{\Delta}^N \arrow{r} \arrow{d} & \mathbb{P}_{T}^N \arrow{d} \\
    0 \arrow{r} &\Delta \arrow{r} & T
\end{tikzcd} \]
where $\Delta$ is the spectrum of the dual numbers Spec($\mathbb{C}[\epsilon]$) whose unique closed point maps to $0 \in T$. (It is clear from the diagram that such a $\Phi$ already satisfies condition (2) of both parts of the theorem.) \par
\noindent\textit{Proof of Step 2.} Let $L = \varphi^{*}(\mathscr{O}_{\mathbb{P}^N} (1))$ and $\widetilde{L} = \widetilde{\varphi}^{*}(\mathscr{O}_{\mathbb{P}_{\Delta}^N}(1))$. \textcolor{blue}{}
Note that $\varphi$ and $\widetilde{\varphi}$ is induced by a sublinear series of $L$ and $\widetilde{L}$ respectively. Then the pair $(\widetilde{X}, \widetilde{L})$ is a first order deformation of the pair $(X,L)$ and hence corresponds to an element $\nu$ of $H^1(\mathscr{E}_L)$. Now the Kuranishi space of deformations, $\textrm{Def}(X,L)$ of $(X,L)$ is smooth (see for example \cite{Hu1}, Section 2.3) and $T_0\textrm{Def}(X,L) \cong H^1(\mathscr{E}_L)$. Hence we can choose a smooth affine algebraic curve $T$ from $\textrm{Def}(X,L)$ passing through the point $(X,L)$ with tangent vector $\nu$. Pulling back the family of deformations of $(X,L)$ on $\textrm{Def}(X,L)$ to $T$, we get a pair ($\mathscr{X}, \mathscr{L}$) consisting of a smooth irreducible family of hyperk\"ahler manifolds $\mathscr{X}$ proper and flat over $T$ and an invertible sheaf $\mathscr{L}$ on $\mathscr{X}$ 
which induces a morphism $\Phi_M : \mathscr{X} \longrightarrow \mathbb{P}_T^{M}$ such that the following squares are Cartesian,\[
\begin{tikzcd}
   X \arrow{r} \arrow[d, "\Phi"] & \widetilde{X} \arrow{r} \arrow[d, "\widetilde{\Phi}"] & \mathscr{X} \arrow[d, "\Phi_M"]\\
    \mathbb{P}^M \arrow{r} \arrow{d} & \mathbb{P}_{\Delta}^M \arrow{r} \arrow{d} & \mathbb{P}_{T}^M \arrow{d} \\
   0 \arrow{r} & \Delta \arrow{r} & T
\end{tikzcd} \]
where $\Phi$ and $\widetilde{\Phi}$ are induced by the complete linear series of $L$ and $\widetilde{L}$ respectively. Now we know that the morphism $\widetilde{\varphi}$ is a composition of $\widetilde{\Phi}$ with a linear projection $\widetilde{\rho}_{\Delta}: \mathbb{P}_{\Delta}^M \rightarrow \mathbb{P}_{\Delta}^n $ (for some $n$) followed by a linear embedding $\widetilde{k}_{\Delta}: \mathbb{P}_{\Delta}^n \hookrightarrow \mathbb{P}_{\Delta}^N $ (since $\widetilde{\varphi}$ is given by a sublinear series of $\widetilde{L}$ followed by a linear embedding). Since 
we can lift the projection $\widetilde{\rho}_{\Delta}$ to a projection $\rho_{T}: \mathbb{P}_{T}^M \rightarrow \mathbb{P}_{T}^n$ and the linear embedding $\widetilde{k}_{\Delta}$ to a linear embedding $k_{T}: \mathbb{P}_{T}^n \hookrightarrow \mathbb{P}_{T}^N$, we have the required Cartesian diagram by composing. \[
\begin{tikzcd}
    X \arrow{r} \arrow[d, "\varphi"] & \widetilde{X} \arrow{r} \arrow[d, "\widetilde{k}_{\Delta} \circ \widetilde{\rho}_{\Delta} \circ \widetilde{\Phi} = \widetilde{\varphi}"]  & \mathscr{X} \arrow[d, "k_{T} \circ \rho_{T} \circ \Phi_M"]\\
   \mathbb{P}^N \arrow{r} \arrow{d} & \mathbb{P}_{\Delta}^N \arrow{r} \arrow{d} & \mathbb{P}_{T}^N \arrow{d} \\
    0 \arrow{r} &\Delta \arrow{r} & T
\end{tikzcd} \]
\noindent\underline{\textit{Step 3.}} In this step we combine the results in the previous steps to prove the theorem. \par
\noindent\textit{Proof of Step 3.} Using the previous steps we can construct a smooth irreducible family $\mathscr{X}$ proper and flat over a smooth pointed irreducible algebraic curve $(T,0)$ and a $T$- morphism $\mathscr{X} \xrightarrow{\Phi} \mathbb{P}_T^N$ such that the fibres $\mathscr{X}_t$ are smooth irreducible projective hyperk\"ahler varieties for all $t\neq 0$, $\mathscr{X}_0 = X$ and $\Phi_0 = \varphi$ and the first order deformation $\widetilde{\varphi}$ of $\varphi$ obtained by pulling back the family to the first infinitesimal neighbourhood $\Delta$ of $0$ in $T$ maps to the class of a homomorphism in $\textrm{Hom}(\mathscr{I}/\mathscr{I}^2, \mathscr{E})$ of rank $ \geq 1$ in part (i) (resp. rank $>\floor{d/2}-1$ in part (i) (b)) and to the class of $\widetilde{Y}$ in part (ii).\par 
\indent Hence by Theorem ~\ref{defs and ropes} $(1)$, we have claims ($1$) and ($2$) of part (i) and by Theorem ~\ref{defs and ropes} $(2)$ we have that $\Phi_t$ is a closed immersion for all $t \neq 0$ which proves claims ($1$) and ($2$) of (ii).\par 
\indent Now suppose $\mathscr{Y} = \textrm{Im}\,\Phi$. Since $\mathscr{X}$ is integral, we have that $\mathscr{Y}$ is integral. This combined with the fact that $T$ is smooth implies that $\mathscr{Y}$ is flat over $T$. Thus, to prove claim (3) in part (ii), it is enough to prove that $\mathscr{Y}_0 = (\textrm{Im}\,\widetilde{\varphi})_0$ which also follows from Theorem ~\ref{defs and ropes} $(2)$ .\QEDB \par 

\vspace{5pt}

Now we prove the same when the universal cover is a Calabi-Yau manifold. The proof in this case will be simpler since now the formally semi-universal deformation space of $\varphi$ is algebraic.  
\begin{theorem}\phantomsection\label{smoothingcy}
Suppose $Y$ be an Enriques manifold of dimension $2n$ and $\pi : X \to Y$ be a universal cover of $Y$ which is a Calabi-Yau manifold. Let $i:Y\hookrightarrow\mathbb{P}^N$ be an embedding of $Y$ into the projective space. Let $ \varphi $ be $i \circ \pi $. 
\begin{itemize}
\item[(i)] There exists a smooth irreducible family $\mathscr{X}$ proper and flat over a smooth pointed affine algebraic curve $(T,0)$ and a $T$-morphism $\Phi:\mathscr{X}\to\mathbb{P}^N_T $ with the following properties:
\begin{itemize}
    \item[(1)] the fibres $\Phi_t: \mathscr{X}_t \to \mathbb{P}^N$, $t \in T-\{0\}$ are finite birational morphisms, 
    \item[(2)] the central fibre $\Phi_0:\mathscr{X}_0 \to \mathbb{P}^N$ is the morphism $\varphi:X \to \mathbb{P}^N$.
\end{itemize}

\item[(ii)] 
Let $\widetilde{Y}$ be an embedded $K$- trivial ribbon inside $\mathbb{P}^N$ supported on $Y$ with conormal bundle $K_Y$ (which exists under the assumptions of Proposition ~\ref{countingcy}, (3)).
Then there exists a smooth irreducible family $\mathscr{X}$ proper and flat over a smooth pointed affine algebraic curve $(T,0)$ and a $T$-morphism $\Phi:\mathscr{X}\to\mathbb{P}^N_T $ with the following properties:
\begin{itemize}
    \item[(1)] the fibres $\Phi_t: \mathscr{X}_t \to \mathbb{P}^N$, $t \in T-\{0\}$ are closed immersions of smooth Calabi-Yau manifolds,
    \item[(2)] the central fibre $\Phi_0:\mathscr{X}_0 \to \mathbb{P}^N$ is the morphism $\varphi:X \to \mathbb{P}^N$.
    \item[(3)] $\textrm{Im}\,\Phi$ is a flat family of schemes over $T$ embedded inside $\mathbb{P}^N_T$ such that for $t \neq 0$, $(\textrm{Im}\,\Phi)_t$ is a smooth Calabi-Yau manifold and the central fibre $(\textrm{Im}\,\Phi)_0$ is $\widetilde{Y}$.
\end{itemize}
\end{itemize}
\end{theorem}
\noindent\textit{Proof.} Using same arguments as in the first step of Theorem ~\ref{smoothinghk} we conclude the existence of first order deformations  $\widetilde{\varphi}: \widetilde{X} \to \mathbb{P}_{\Delta}^N$ of $\varphi: X \to \mathbb{P}^N$ such that their class $\widetilde{\varphi} \in H^0(\mathscr{N}_{\varphi})$ maps to the class of a homomorphism in $\textrm{Hom}(\mathscr{I}/\mathscr{I}^2, \mathscr{E})$ of rank $ \geq 1$ in part (i) and to the class of  $\widetilde{Y} \hookrightarrow \mathbb{P}^N$ in  $\textrm{Hom}(\mathscr{I}/\mathscr{I}^2,K_Y)$ in part (ii).\par 
\indent Now since $H^2(\mathscr{O}_X) = 0$ we conclude that $X$ has an algebraic formally semi-universal deformation space. Again since $\varphi$ is non-degenerate we conclude that $\varphi$ has an algebraic formally semi-universal deformation as well (see \cite{GGP1}, Remark $1.7$). Now by Lemma ~\ref{unobphicy}, $\varphi$ is unobstructed. So we can choose a smooth algebraic curve $T$ from the smooth semi-universal deformation space of $\varphi$ with tangent vector $\widetilde{\varphi}$. Pulling back the semi-universal family along $T$ we get a deformation $\mathscr{X} \to \mathbb{P}^N_T$ of $\varphi$. Hence by Theorem ~\ref{defs and ropes} $(1)$, we have claims ($1$) and ($2$) of part (i) and by Theorem ~\ref{defs and ropes} $(2)$ we have that claims ($1$), ($2$) and ($3$) of part (ii).\QEDB\par 

\vspace{5pt}

Now we prove a consequence of smoothing, namely the smoothness of embedded (regular) $K$- trivial ribbons on an Enriques manifold with conormal bundle $\mathscr{E}=K_Y$ in its Hilbert scheme.
\begin{proposition}\phantomsection\label{smhs}
Let $\widetilde{Y}\hookrightarrow\mathbb{P}^N$ be a projective $K$- trivial ribbon on an Enriques manifold $Y$ of dimension $2n$ and index $d=2$ extending an  embedding $Y\hookrightarrow\mathbb{P}^N$. 
Then the Hilbert point of $\widetilde{Y}$ inside $\mathbb{P}^N$ is nonsingular.
\end{proposition}
\noindent\textit{Proof.} Here we give the proof when the universal cover $X$ is a hyperk\"ahler manifold. The proof of the case when $X$ is Calabi-Yau is similar. Suppose the sectional genus of $Y\hookrightarrow \mathbb{P}^N$ is $g$. Since $\widetilde{Y}$ admits an embedded smoothing, first we compute the dimension of the component that parametrizes smooth hyperk\"ahler manifolds inside $\mathbb{P}^N$. So our situation is as follows: we have a hyperk\"ahler manifold $X'\hookrightarrow\mathbb{P}^N$ such that $h^0(\mathscr{O}_{X'}(1))=2g$. We will calculate $h^0(\mathscr{N}_{X'/\mathbb{P}^N})$.
It is easy to check that we have the following
\begin{align}\phantomsection\label{dimhshk}
    h^0(\mathscr{N}_{X'/\mathbb{P}^N})=2g(N+1)+h^{1,1}_{X'}-2.
\end{align} 
Indeed, to get ~\eqref{dimhshk} one has to use the following exact sequence
\begin{align}
    0\to T_{X'}\to T_{\mathbb{P}^N}\vert_{X'}\to \mathscr{N}_{X'/\mathbb{P}^N}\to 0
\end{align}
and the restriction of the Euler sequence of $\mathbb{P}^N$ on $X'$ keeping in mind that the connecting homomorphism  $H^1(\mathscr{N}_{X'/\mathbb{P}^N})\to H^2(T_{X'})$ is an isomorphism (see Lemma ~\ref{unobofphi}).\par 
\indent Therefore the Hilbert point of $\widetilde{Y}$ will be nonsingular if and only if $h^0(\mathscr{N}_{\widetilde{Y}/\mathbb{P}^N})=2g(N+1)+h^{1,1}_X-2$.

\indent Next, we use the following exact sequences (see \cite{GGP}, Lemma 4.2) to calculate $h^0(\mathscr{N}_{\widetilde{Y}/\mathbb{P}^N})$
\begin{align}\phantomsection\label{adh1}
    0\to \mathscr{N}_{\widetilde{Y}/\mathbb{P}^N}\vert_Y\otimes K_Y\to \mathscr{N}_{\widetilde{Y}/\mathbb{P}^N}\to\mathscr{N}_{\widetilde{Y}/\mathbb{P}^N}\vert_Y\to 0,
\end{align}
\begin{align}\phantomsection\label{adh2}
     0\to K_Y\to\mathscr{N}_{Y/\mathbb{P}^N}\to\mathscr{H}om_Y(\mathscr{I}_{\widetilde{Y}/\mathbb{P}^N}/\mathscr{I}^2_{Y/\mathbb{P}^N},\mathscr{O}_Y)\to 0,
\end{align}
\begin{align}\phantomsection\label{adh3}
    0\to \mathscr{H}om_Y(\mathscr{I}_{\widetilde{Y}/\mathbb{P}^N}/\mathscr{I}^2_{Y/\mathbb{P}^N},\mathscr{O}_Y)\to \mathscr{N}_{\widetilde{Y}/\mathbb{P}^N}\vert_Y\to\mathscr{O}_Y\to 0.
\end{align}
The exact sequence ~\eqref{adh1} gives us the following inequality
\begin{align}\phantomsection\label{smeq1}
    h^0(\mathscr{N}_{\widetilde{Y}/\mathbb{P}^N})\leq h^0(\mathscr{N}_{\widetilde{Y}/\mathbb{P}^N}\vert_Y\otimes K_Y)+h^0(\mathscr{N}_{\widetilde{Y}/\mathbb{P}^N}\vert_Y).
\end{align}
We use ~\eqref{adh2} and ~\eqref{adh3} to see that
\begin{align}\phantomsection\label{smeq2}
    h^0(\mathscr{N}_{\widetilde{Y}/\mathbb{P}^N}\vert_Y\otimes K_Y)=h^0(\mathscr{N}_{Y/\mathbb{P}^N}\otimes K_Y)-1=g(N+1)+h^1(T_Y\otimes K_Y)-2
\end{align} 
\begin{align}\phantomsection\label{smeq3}
    h^0(\mathscr{N}_{\widetilde{Y}/\mathbb{P}^N}\vert_Y)=g(N+1)+h^1(T_Y)-1\quad \textrm{or} \quad g(N+1)+h^1(T_Y).
\end{align} 
Using ~\eqref{smeq1}, ~\eqref{smeq2} and ~\eqref{smeq3} we deduce the following inequality
\begin{align}
    h^0(\mathscr{N}_{\widetilde{Y}/\mathbb{P}^N})\leq 2g(N+1)+h^1(T_X)-2.
\end{align} 
Our assertion follows since $h^0(\mathscr{N}_{\widetilde{Y}/\mathbb{P}^N})$ is at least the dimension of the component.\QEDB

\section{Appendix}\label{App} 

Here we discuss the condition that an embedding $Y\hookrightarrow\mathbb{P}^N$ of the Enriques manifold satisfies $N\geq 2\dim(Y)+1$. It is known that any variety $Y$ can be embedded inside $\mathbb{P}^{2\dim(Y)+1}$ (see \cite{Sha}). We show that $N\geq 2\dim(Y)+1$ is always the case if $Y$ is an Enriques manifold of index two and $X$ is one of the known classes of hyperk\"ahler manifolds of dimension six . 

\begin{proposition}\phantomsection\label{optwh}
Let $Y$ be an Enriques manifold of dimension six and index two whose universal double cover $X$ is either $K3^{[3]}$ or $K^3(A)$. Let $Y\hookrightarrow \mathbb{P}^N$ be an embedding. Then $N\geq 13$.
\end{proposition}
\noindent\textit{Proof.}
It is enough to show that $Y$ can not be embedded inside $\mathbb{P}^{12}$. For the sake of contradiction, assume $Y\hookrightarrow\mathbb{P}^{12}$ is a degree $d$ embedding. We have  that $\pi:X\to Y$ is the universal cover. Let $\mathscr{N}$ be the normal bundle of the embedding $Y\hookrightarrow\mathbb{P}^{12}$ and $K$ be the canonical bundle of $Y$.\par 
\indent From the Euler sequence on $\mathbb{P}^{12}$ we get that 
\begin{align}
    c_t(T_{\mathbb{P}^{12}})=(1+ht)^{13}
\end{align}
where $h\in A^1(\mathbb{P}^{12})$ is the class of a hyperplane. \par 
\indent Suppose $\mathscr{N}$ is the normal bundle for the embedding. From the following exact sequence 
\begin{align}
    0\to T_Y\to T_{\mathbb{P}^{12}}\vert_Y\to\mathscr{N}\to 0
\end{align}
we see that $c_t(T_Y)\cdot c_t(\mathscr{N})=c_t(T_{\mathbb{P}^{12}}\vert_Y)$. That gives us 
\begin{align}
    \left(1-Kt+\displaystyle\sum_{i=2}^{i=6}c_i(Y)t^i\right)\left(1+\displaystyle\sum_{i=1}^{i=6}c_i(\mathscr{N})t^i\right)=\displaystyle\sum_{i=0}^{i=6}\binom{13}{i}H^it^{6-i}
\end{align}
where $H\in A^1(Y)$ is the class of a hyperplane section.\par 
\indent Expanding the equality above we get the following equations
\begin{align}\phantomsection\label{embeq1}
    c_1(\mathscr{N})=13H+K
\end{align}
\begin{align}\phantomsection\label{embeq2}
    c_2(\mathscr{N})=\binom{13}{2}H^2+Kc_1(\mathscr{N})-c_2(Y)
\end{align}
\begin{align}\phantomsection\label{embeq3}
    c_3(\mathscr{N})=\binom{13}{3}H^3+Kc_2(\mathscr{N})-c_2(Y)c_1(\mathscr{N})-c_3(Y)
\end{align}
\begin{align}\phantomsection\label{embeq4}
    c_4(\mathscr{N})=\binom{13}{4}H^4+Kc_3(\mathscr{N})-c_2(Y)c_2(\mathscr{N})-c_3(Y)c_1(\mathscr{N})-c_4(Y)
\end{align}
\begin{align}\phantomsection\label{embeq5}
    c_5(\mathscr{N})=\binom{13}{5}H^5+Kc_4(\mathscr{N})-c_2(Y)c_3(\mathscr{N})-c_3(Y)c_2(\mathscr{N})-c_4(Y)c_1(\mathscr{N})-c_5(Y)
\end{align}
\begin{align}\phantomsection\label{embeq6}
\begin{split}
    \binom{13}{6}H^6=c_6(\mathscr{N})-Kc_5(\mathscr{N})+c_2(Y)c_4(\mathscr{N})+c_3(Y)c_3(\mathscr{N})\\+c_4(Y)c_2(\mathscr{N})+c_5(Y)c_1(\mathscr{N})+c_6(Y)
    \end{split}
\end{align}
Next we put the values of $c_i(\mathscr{N})$ for $i=1,\dots,5$ in ~\eqref{embeq6} and pull them back to $X$ to get the following (recall that the odd Chern classes of $X$ are trivial and $c_6(\mathscr{N})=Y.Y=d^2$)
\begin{align}\phantomsection\label{embeq7}
\begin{split}
    \binom{13}{6}\dfrac{(\pi^*H)^6}{2}=\left(\dfrac{(\pi^*H)^6}{2}\right)^2+\binom{13}{4}\dfrac{(\pi^*H)^4}{2}c_2(X)+\binom{13}{2}\dfrac{(\pi^*H)^2}{2}c_4(X)-\binom{13}{2}\dfrac{(\pi^*H)^2}{2}\left(c_2(X)\right)^2\\+\dfrac{1}{2}\left(\left(c_2(X)\right)^3-2c_2(X)c_4(X)+c_6(X)\right).
\end{split}    
\end{align}
Let $q(-,-)$ and $c$ be the Beauville-Bogomolov-Fujiki form and Fujiki constant respectively on $X$. We know the value of the constants $a_i$ (see  ~\eqref{SRR}) if $X$ is of type $K3^{[3]}$ or $K^3(A)$. \par 
\indent Suppose $X$ is of the form $K3^{[3]}$. The following Riemann-Roch formula has been proven by  Ellingsrud, G\"ottsche and Lehn (see \cite{Ell}) when $X$ is of $K3^{[3]}$ type. Fujiki constant takes the value 15.
\begin{align}\phantomsection\label{rrk33}
    \chi(L) = \binom{\frac{1}{2}q(L)+4}{3}.
\end{align}
Comparing degree-two terms of the above expression and the usual Riemann-Roch  we find that
\begin{align}
    c_2(X)L^4=108q(L).
\end{align}
Note that this gives us the value of the constant $C(c_2)$ as well (see Theorem ~\ref{GRR} for the definition).
\begin{align}\phantomsection\label{cc2}
    C(c_2)=108
\end{align}
Comparing degree-one terms of ~\eqref{rrk33} and the usual Riemann-Roch expression we get
\begin{align}\phantomsection\label{embeq8}
    3(c_2(X))^2L^2-c_4(X)L^2=3120q(L).
\end{align}
Recall that we have the following equality on hyperk\"ahler manifolds of dimension $2n$ for a line bundle $L$ (see \cite{Nie} and \cite{CaoJ}, proof of Theorem 3.2)
\begin{align}\phantomsection\label{embeq9}
    \int_X\sqrt{\textrm{td}(X)}L^{2n-4}=(2n-4)!\binom{n}{n-2}\lambda(L)^{n-2}\int_X\sqrt{\textrm{td}(X)} 
\end{align}
where $\lambda(L)$ is the characteristic value given by the following expression
\begin{align}\phantomsection\label{lambda}
    \lambda(L)=\dfrac{12c}{(2n-1)C(c_2(X))}q(L).
\end{align}
Comparing degree-one terms of the left hand side and right hand side of ~\eqref{embeq9} and using ~\eqref{cc2} and the fact that $\int_X\sqrt{\textrm{td}(X)}=\dfrac{9}{16}$ (see \cite{Sa}, Proposition 19) we get the following equality
\begin{align}\phantomsection\label{embeq10}
    7(c_2(X))^2L^2-4c_4(X)L^2=6480q(L).
\end{align}
We solve ~\eqref{embeq8} and ~\eqref{embeq10} to get that
$(c_2(X))^2L^2=1200q(L)\textrm{ and }c_4(X)L^2=480q(L)$.\par 
\indent Suppose $X$ is of the form $K^3(A)$. The Riemann-Roch formula for this case is given below, it was calculated by Britze-Nieper (see \cite{Br}). The value of the Fujiki constant in this case is 60.
\begin{align}
    \chi(L) = 4 \binom{\frac{1}{2}q(L)+3}{3}.
\end{align}
We carry out the same computations ($\int_X\sqrt{\textrm{td}(X)}=\dfrac{2}{3}$ by \cite{Sa}, Proposition 21) to get that 
\begin{align}
    (c_2(X))^2L^2=1920q(L)\textrm{ and }c_4(X)L^2=480q(L).
\end{align}
Now, $c_2(X)^3-2c_2(X)c_4(X)+c_6(X)$ can be calculated form the Hodge diamond (see \cite{Sa}, Appendix B). The Hodge diamonds of $K3^{[3]}$ and $K^3(A)$ are given below (see \cite{Bel}).
\begin{multicols}{2}
\tiny{
\begin{center}
    Hodge diamond of $K3^{[3]}$
\end{center}\[
\begin{array}{ccccccccccccc}
   & & & & & & 1 & & & & & &  \\
   & & & & & 0 & & 0 & & & & & \\
   & & & & 1 & & 21 & & 1 & & & & \\
   & & & 0 & & 0 & & 0 & & 0 & & & \\
   & & 1 & & 22 & & 253 & & 22 & & 1 & & \\
   & 0 & & 0 & & 0 & & 0 & & 0 & & 0 & \\
   1 & & 21 & & 253 & & 2004 & & 253 & & 21 & & 1 \\
   & 0 & & 0 & & 0 & & 0 & & 0 & & 0 & \\
   & & 1 & & 22 & & 253 & & 22 & & 1 & & \\
   & & & 0 & & 0 & & 0 & & 0 & & & \\
   & & & & 1 & & 21 & & 1 & & & & \\
   & & & & & 0 & & 0 & & & & & \\
   & & & & & & 1 & & & & & &  \\
\end{array}\]}
\columnbreak
\begin{center}
    Hodge diamond of $K^3(A)$
\end{center}
\tiny{\[
\begin{array}{ccccccccccccc}
   & & & & & & 1 & & & & & &  \\
   & & & & & 0 & & 0 & & & & & \\
   & & & & 1 & & 5 & & 1 & & & & \\
   & & & 0 & & 4 & & 4 & & 0 & & & \\
   & & 1 & & 6 & & 37 & & 6 & & 1 & & \\
   & 0 & & 4 & & 24 & & 24 & & 4 & & 0 & \\
   1 & & 5 & & 37 & & 372 & & 37 & & 5 & & 1 \\
   & 0 & & 4 & & 24 & & 24 & & 4 & & 0 & \\
   & & 1 & & 6 & & 37 & & 6 & & 1 & & \\
   & & & 0 & & 4 & & 4 & & 0 & & & \\
   & & & & 1 & & 5 & & 1 & & & & \\
   & & & & & 0 & & 0 & & & & & \\
   & & & & & & 1 & & & & & &  \\
\end{array}\]}
\end{multicols}
\indent Let $\chi^p(X)=\sum (-1)^qh^{p,q}(X)$. Then \cite{Sa}, Appendix B gives the following equations.
\begin{align}
    c_2(X)^3=7272\chi^0(X)-184\chi^1(X)-8\chi^2(X),
\end{align}
\begin{align}
    c_2(X)c_4(X)=1368\chi^0(X)-208\chi^1(X)-8\chi^2(X),
\end{align}
\begin{align}
    c_6(X)=36\chi^0(X)-16\chi^1(X)+4\chi^2(X).
\end{align}
We set $L=\pi^*H$ and $x=q(\pi^*H)$. We put the values of $(c_2(X))^2$, $c_4(X)L^2$ and $c_2(X)^3-2c_2(X)c_4(X)+c_6(X)$ calculated above in ~\eqref{embeq6} to get an equation in $x$.\par 
\indent If $X$ is of the form $K3^{[3]}$ the equation we obtain is the following
\begin{align}\phantomsection\label{embeq13}
    51480x^3=225x^6+154440x^2-112320x+21120.
\end{align}
\indent When $X$ is of the form $K^3(A)$ the equation we obtain is given below
\begin{align}\phantomsection\label{embeq14}
    51480x^3=900x^6+102960x^2-56160x+8664.
\end{align}
That concludes the proof since neither ~\eqref{embeq13} nor ~\eqref{embeq14} has positive even integer solution.\QEDB\par 

\vspace{5pt}

Now we proceed to prove the same when the universal cover $X$ is O'Grady's six dimensional hyperk\"ahler manifold $\mathscr{M}_6$. In order to do that we follow the same procedure but first we need $\int_{\mathscr{M}_6}\sqrt{\textrm{td}(\mathscr{M}_6)}$ which we calculate below.\par 
\begin{lemma}
$\int_{\mathscr{M}_6}\sqrt{\textrm{td}(\mathscr{M}_6)}=\dfrac{2}{3}$.
\end{lemma}
\noindent\textit{Proof.} We know the Hodge diamond of $\mathscr{M}_6$ by  Mongardi-Rapagnetta-Sacc\`a’s computation (see \cite{Mon}).
\begin{center}
    \tiny{Hodge diamond of $\mathscr{M}_6$
\[
\begin{array}{ccccccccccccc}
   & & & & & & 1 & & & & & &  \\
   & & & & & 0 & & 0 & & & & & \\
   & & & & 1 & & 6 & & 1 & & & & \\
   & & & 0 & & 0 & & 0 & & 0 & & & \\
   & & 1 & & 12 & & 173 & & 12 & & 1 & & \\
   & 0 & & 0 & & 0 & & 0 & & 0 & & 0 & \\
   1 & & 6 & & 173 & & 1144 & & 173 & & 6 & & 1 \\
   & 0 & & 0 & & 0 & & 0 & & 0 & & 0 & \\
   & & 1 & & 12 & & 173 & & 12 & & 1 & & \\
   & & & 0 & & 0 & & 0 & & 0 & & & \\
   & & & & 1 & & 6 & & 1 & & & & \\
   & & & & & 0 & & 0 & & & & & \\
   & & & & & & 1 & & & & & &  \\
\end{array}\]}
\end{center}
As before $\chi^p(X)=\sum (-1)^qh^{p,q}(X)$. Then we have $\chi^0=4$, $\chi^1=-24$, $\chi^2=348$. \par 
\indent It has been proven by Sawon (see \cite{Sa}) that
$\int_{\mathscr{M}_6}\sqrt{\textrm{td}(\mathscr{M}_6)}=-\dfrac{1}{48^3\cdot3!}\left(s_2^3+\dfrac{12}{5}s_2s_4+\dfrac{64}{35}s_6\right)$ where
\begin{align*}
    s_2^3=-58176\chi^0+1472\chi^1+64\chi^2,\quad s_2s_4=-18144\chi^0-928\chi^1-32\chi^2,
\end{align*}
\begin{align*}
    s_6=-6552\chi^0-784\chi^1-56\chi^2.
\end{align*}
A direct computation yields the result.\QEDB\par 

\vspace{5pt}

Once we know the value of $\int_{\mathscr{M}_6}\sqrt{\textrm{td}(\mathscr{M}_6)}$ we can carry out the procedure explained in the proof of Proposition ~\ref{optwh}.
\begin{proposition}\phantomsection\label{optwhm6}
Let $Y$ be an Enriques manifold and $Y\hookrightarrow\mathbb{P}^N$ be a closed embedding. Assume that the universal cover $X$ of $Y$ is O'Grady's six dimensional hyperk\"ahler manifold $\mathscr{M}_6$. Then $N\geq 13$.
\end{proposition}
\noindent\textit{Proof.} Recall that on $\mathscr{M}_6$ the Riemann-Roch expression is the following (see \cite{CaoJ})
\begin{align}\phantomsection\label{RROG}
\chi(L)=4\binom{\lambda(L)+3}{3}
\end{align}
where $\lambda(L)$ is the characteristic value given by ~\eqref{lambda}.\par 
\indent Comparing degree-one terms of this expression and the usual Riemann-Roch expression, we get 
\begin{align}
    \dfrac{1}{12}c_2(\mathscr{M}_6)\cdot\dfrac{1}{24}L^4=4(\lambda(L))^2
\end{align}
Using $c_2(\mathscr{M}_6)L^4=C(c_2(\mathscr{M}_6))q(L)$ and the value of the Fujiki constant (which is 60), we get that $C(c_2(\mathscr{M}_6))=288$. We use them to get the following two equations; the first one is obtained by comparing degree-two terms of ~\eqref{RROG} and the usual Riemann-Roch expression, the second one is a consequence of ~\eqref{embeq9}.
\begin{align}
    3(c_2(\mathscr{M}_6))^2L^2-c_4(\mathscr{M}_6)L^2=5280q(L) 
\end{align}
\begin{align}
    7(c_2(\mathscr{M}_6))^2L^2-4c_4(\mathscr{M}_6)L^2=11520q(L)
\end{align}
We solve the above two equations to get $(c_2(\mathscr{M}_6))^2L^2=1920q(L)$ and $c_4(\mathscr{M}_6)L^2=480q(L)$. Moreover the values of $(c_2(\mathscr{M}_6))^3$, $c_2(\mathscr{M}_6)c_4(\mathscr{M}_6)$ and $c_6(\mathscr{M}_6)$ have been calculated in \cite{Mon}, Corollary 6.8. We put these values in ~\eqref{embeq7} by setting $L=\pi^*(H)$ where $H$ is the class of a hyperplane section of $Y$ and $\pi$ is the map from the universal cover. As before, setting $x=q(\pi^*(H))$, we obtain 
\begin{align}
    51480x^3=900x^6+102960x^2-56160x+8640
\end{align}
which has no positive integer solution. That concludes the proof.\QEDB\par

\bibliographystyle{plain}

\end{document}